\documentclass[11pt]{amsart} % change to \documentclass[a4paper,11pt]{amsart} if necessary. The current one is suitable for legal size paper.
\usepackage[lmargin=1in,rmargin=1in,tmargin=1in,bmargin=1in]{geometry}
\usepackage{tikz}
\usepackage{graphicx, float, epstopdf}
\usepackage{hyperref, color}
\usepackage{comment}
\usepackage{centernot}
\usepackage{fancyhdr}
\usepackage[utf8]{inputenc}

\usepackage{amsfonts,amssymb,amsmath,amsthm,mathrsfs}
\usepackage{graphics, setspace}
\usepackage{braket}
\usepackage{mathtools}
\usepackage{longtable}
\usepackage{pythonhighlight}
\usepackage{datetime}
\usepackage{algorithmic,algorithm}
\numberwithin{equation}{section}
\numberwithin{figure}{section}
\allowdisplaybreaks[4]                        %allows breaking multiline environments in amsmath commands
 
\newtheorem{lemma}{Lemma}[section]
\newtheorem{theorem}{Theorem}[section]

\newtheorem{example}{Example}[section]

\newtheorem{corollary}[lemma]{Corollary}
\theoremstyle{definition}
\newtheorem{definition}{Definition}[section]
\newtheorem{remark}{Remark}[section]
\newcommand{\Z}{\mathbb{Z}}

\newcommand{\R}{\mathbb{R}}
\newcommand{\C}{\mathbb{C}}
\newcommand{\N}{\mathbb{N}}

\newcommand{\e}{\operatorname{e}}

\newcommand{\real}{\operatorname{Re}}

\newcommand{\modu}{\operatorname{mod}}

\newcommand{\Pri}{\mathbb{P}}

\usepackage{tikz}

\begin{document}

\title{exponential sums weighted by additive functions}

\author[A.~Gafni]{Ayla Gafni}
\address{Department of Mathematics, The University of Mississippi, P.O. Box 1848, University, MS 38677}
\email{argafni@olemiss.edu}
\author[N.~Robles]{Nicolas Robles}
\address{RAND Corporation, Engineering and Applied Sciences, 1200 S Hayes St, Arlington, VA, 22202}
\email{nrobles@rand.org}

\begin{abstract}  
    We introduce a general class $\mathcal{F}_0$ of additive functions $f$ such that $f(p) = 1$ and prove a tight bound for exponential sums of the form $\sum_{n \le x} f(n) \e(\alpha n)$ where $f \in \mathcal{F}_0$ and $\e(\theta) = \exp(2\pi i \theta)$. Both $\omega$, the number of distinct primes of $n$, and $\Omega$, the total number primes of $n$, are members of $\mathcal{F}_0$. As an application of the exponential sum result, we use the Hardy-Littlewood circle method to find the asymptotics of the Goldbach-Vinogradov ternary problem associated to $\Omega$, namely we show the behavior of $r_\Omega(N) = \sum_{n_1+n_2+n_3=N}\Omega(n_1)\Omega(n_2)\Omega(n_3)$, as $N \to \infty$.  Lastly, we end with a discussion of further applications of the main result.
\end{abstract}


\maketitle
%%%%%%%%%%%%%%%%%%%%%%%%%%%%%%%%%%%%%%%%%%%%%%%%%%%%%%%%
%\tableofcontents
\section{Introduction} \label{sec:introduction}

Let $f: \mathbb N \rightarrow \mathbb C$ be an arithmetic function and consider the exponential sum 
\begin{equation}\label{eq:exp sum f}
S_f(\alpha; X) := \sum_{n\le X} f(n) \e(\alpha n),
\end{equation} 
for $\alpha \in \mathbb R$ and $X$ a positive real number.  Here we use the standard notation $\e(\theta) = e^{2\pi i \theta}$. Exponential sums of this form are ubiquitous throughout mathematics, and obtaining nontrivial bounds for these sums is central to many topics in analytic number theory.  Indeed, at the heart of the Hardy-Littlewood circle method is Weyl's inequality, which can be viewed as a bound on $S_f(\alpha; X)$ where $f$ is the indicator function of a monic polynomial \cite[Theorem 4.1]{Vaughan1997}.  

If $f=\Lambda$, the von Mangoldt function,  Vinogradov \cite{Vi42} showed that
\[
S_{\Lambda}(\alpha, X) \ll \bigg(\frac{X}{q^{1/2}}+X^{4/5}+X^{1/2}q^{1/2}\bigg)(\log X)^4,
\]
for real $\alpha$ such that $|\alpha-\frac{a}{q}|\le q^{-2}$ with $(a,q)=1$. In the same spirit, the case $f=\mu$, the M\"{o}bius function, was considered by Davenport \cite{davenportInequality}, and $f=\mu*\mu$ was considered in \cite{brz23}. The situation for $f=\sigma_{k_1,k_2}$, the generalized divisor function, can be found in \cite{MotohashiDivisor} for $k_1=k_2=0$ and in \cite{BRZZ23}. General bounds when $f$ is a multiplicative function with certain growth conditions were treated in \cite{MV77}.

Our objective in this note is to study \eqref{eq:exp sum f} when $f$ is an additive function and solve an associated problem in the theory of partitions as an application.

%%%%%%%%%%%%%%%%%%%%%%%%%%%
\subsection{Exponential sums with additive coefficients}
One of the main results of this paper is a bound on $S_f(\alpha; X)$ for $f$ belonging to a broad class of additive functions. Recall that an arithmetic function $f$ is called additive if it satisfies the condition $f(ab)=f(a)+f(b)$ when $(a,b)=1$.  

\begin{definition} \label{def: class F}
We define $\mathcal F_0$ to be the set of all additive functions $f$ satisfying
$f(p)= 1$ for each prime $p$.
\end{definition}

\begin{example} \label{example: Omega k}
    Let $n=p_1^{\alpha_1}p_2^{\alpha_2}\cdots p_r^{\alpha_r}$ denote the representation of $n$ as a product of powers of distinct primes and define $\Omega_k(n)=\alpha_1^k+\alpha_2^k+\cdots+\alpha_r^k$ for fixed $k\ge 0$, see \cite{duncanClass1962}. The function $f=\Omega_k$ is a  significant example of $f$ in the class $\mathcal{F}_0$. In particular, $\Omega_0=\omega$, the number of distinct prime divisors of $n$ and $\Omega_1=\Omega$, the total number of prime divisors of $n$, honoring multiplicity.
\end{example}

\begin{remark} \label{rem:convolution_structure}
We denote the set of primes by $\mathbb{P}$ and moreover we write $\mathbb{P}^* = \{p^k : p \in \mathbb{P},\, k \in \N^{\ge 1}\}$ for the set of prime powers.
An arithmetic function $f$ is additive if and only if $f = h * 1$
for some function $h$ supported on prime powers. The condition
$f \in \mathcal{F}_0$ then amounts to $h(p) = 1$ for all primes $p$,
with the values $\{h(p^k)\}_{k \ge 2}$ left unconstrained. For
instance, setting $h = \mathbf{1}_{\mathbb{P}}$ gives $\omega$ whereas setting
$h(p^k) = 1$ for all $k \ge 1$ yields $\Omega$.
\end{remark}

The class $\mathcal{F}_0$ is related, but not identical, to the idea of `small' and `large' additive functions. `Small' additive functions, as described in \cite[$\mathsection$ 14.7]{IvicBook}, are typically functions like $\omega(n)=\sum_{p|n}1$ and $\Omega(n) = \sum_{p^a \| n}a$, whose average order is $C \log \log n$. In contrast, `large' additive functions are typically functions such as
\begin{align}
    \beta(n) = \sum_{p|n} p \quad \textnormal{and} \quad B(n) = \sum_{p^a\|n} a p.
\end{align}
For further discussions, see \cite{deKoninckIvic, kubiliusBook}. 

 Our result on exponential sums of the form \eqref{eq:exp sum f} with $f \in \mathcal{F}_0$ is as follows.
\begin{theorem} \label{thm:exp sum bound class F}
    Let $\alpha \in \mathbb R$, $a \in \mathbb Z$ and $q \in \mathbb N$ such that $(a,q)=1$ and 
$\left|\alpha-a/q\right| \le q^{-2}$.  Let $f(n)$ be an additive function in $\mathcal F_0$, and let $F_f(X)$ denote the least upper bound of $f$ on powers of primes up to $X$. That is,
\begin{align}
    F_f(X) = \max_{p^\ell\le X} |f(p^\ell)|.
\end{align} Then for any $\Delta \in (0,\frac{1}{2})$ we have
    \begin{align}
        S_f(\alpha;X) \ll \bigg(\frac{X}{q^{\Delta}} + X^{\frac{5}{6}} + X^{1-\Delta} q^{\Delta}\bigg)\left((\log X)^4 + (\log X)F_f(X)\right).
    \end{align}
\end{theorem}

This result compares well with other similar bounds in the literature. For instance, if we specialize to the representative case $f=\omega$, then Theorem \ref{thm:exp sum bound class F} yields 
    \begin{align} \label{eq: exp sum comparison 1}
        S_\omega(\alpha;X) \ll \bigg(\frac{X}{q^{\Delta}} + X^{\frac{5}{6}} + X^{1-\Delta} q^{\Delta}\bigg) (\log X)^4 
    \end{align}
for any $\Delta \in (0,\frac{1}{2})$. Since $\omega = \mathbf{1} * \mathbf{1}_{\mathbb{P}}$, the techniques put forward in \cite{semiprimes} are very much applicable. While those techniques are somewhat more versatile as they can be applied to a wider range of arithmetic functions, they yield the slightly weaker bound
\begin{align} \label{eq: exp sum comparison 2}
    S_\omega(\alpha;X) \ll \frac{X}{q^{\frac{1}{6}}} (\log X)^{\frac{7}{3}} + X^{\frac{16}{17}} (\log X)^{\frac{39}{17}} + X^{\frac{7}{8}}q^{\frac{1}{8}} (\log X)^{\frac{9}{4}}. 
\end{align}
In comparison, the improvements obtained in \cite{DongRoblesZaharescuZeindler} yield the tighter bound
\begin{align} \label{eq: exp sum comparison 3}
    S_\omega(\alpha;X) \ll \frac{X}{q^{\frac{1}{4}}} (\log X)^{\frac{5}{2}} + X^{\frac{6}{7}} (\log X)^{\frac{19}{7}} + X^{\frac{3}{4}}q^{\frac{1}{4}} (\log X)^{\frac{5}{2}}. 
\end{align}
Moreover, in \cite{madhudas} the methods from \cite{MV77} were adapted from multiplicative functions to additive functions with specific growth conditions. 
For $\omega$ the bound in \cite{madhudas} is
\begin{align} \label{eq: exp sum comparison 4}
    S_\omega(\alpha;X) \ll \bigg(\frac{X}{\log X} + \frac{X (\log R)^{\frac{3}{2}}}{R^{\frac{1}{2}}}\bigg)\log \log X \quad \textnormal{for} \quad 2 \le R \le q \le \frac{X}{R}. 
\end{align}
We note that $\frac{X}{\log X} + \frac{X (\log R)^{3/2}}{R^{1/2}}$ corresponds to the bound for multiplicative functions from \cite{MV77}.

The reader is also referred to \cite[$\mathsection$ 1.2]{DongRoblesZaharescuZeindler} to see how additive functions can be manufactured from other fundamental arithmetic functions and how their associated exponential sums can be bounded. 

For many applications of the circle method, including the application given in this paper, any of the bounds in \eqref{eq: exp sum comparison 1}, \eqref{eq: exp sum comparison 2}, \eqref{eq: exp sum comparison 3} or \eqref{eq: exp sum comparison 4} is enough to handle the minor arcs associated to $\omega$. 
The reason for this is that when setting up the arcs in our applications we only need $q$ to be a power of $\log X$ and therefore minimal savings suffice in this case. However, in general, it is much more advantageous, and difficult, to decrease the exponent associated to powers of $X$ than it is to decrease the exponent associated to powers of $\log X$. This could be particularly useful in future number theoretic applications outside the realm of the circle method.  In addition, Theorem \ref{thm:exp sum bound class F} offers the added benefit of providing a flexible range for $\Delta$ which could be customized for different applications.

\subsection{An application to the circle method} 
To illustrate an application of Theorem \ref{thm:exp sum bound class F}, we introduce the Goldbach-Vinogradov problem.  Vinogradov \cite{Vi42} showed that
\begin{align} \label{eq: Vinogradov}
    S_{\mathbf{1}_{\mathbb{P}}}(\alpha;X)= \sum_{p \le X} \e(\alpha p) \ll \bigg(\frac{X}{\sqrt{q}}+X^{\frac{4}{5}}+\sqrt{X}\sqrt{q}\bigg)(\log X)^3,
\end{align}
where the sum is taken over primes $p$.
The proof was later simplified considerably by Vaughan \cite{Va77}, see also \cite[$\mathsection$ 25]{Da74} and \cite[$\mathsection$ 23]{primesDimitris}. Let $\theta(n) = \mathbf{1}_{\mathbb{P}}(n)\log n$. Employing \eqref{eq: Vinogradov} along with the Hardy-Littlewood circle method, Vinogradov proved that
\begin{align} \label{eq: Vinogradov theta}
r_\theta(N):=\sum_{n_1+n_2+n_3=N} \theta(n_1)\theta(n_2)\theta(n_3) = \mathfrak{S}(N)\frac{N^2}{2}+O\bigg(\frac{N}{(\log N)^A}\bigg),
\end{align}
for every $A>0$ as $N \to \infty$ and where 
\[
\mathfrak{S}(N)=\sum_{q=1}^\infty \frac{\mu(q)c_q(N)}{\varphi(q)^3}. 
\]
Here $c_q(n)$ is the Ramanujan sum and $\varphi(q)$ is the Euler totient function. Using an elementary argument the asymptotic in \eqref{eq: Vinogradov theta} can effectively be rearranged to yield the celebrated result that every sufficiently large odd integer is the sum of three primes.

The object of study in our application of Theorem \ref{thm:exp sum bound class F} is the Goldbach-Vinogradov problem associated to $\Omega$, which we recall is an additive function in $\mathcal{F}_0$.  Following the analogy with \eqref{eq: Vinogradov theta}, let $r_{\Omega}(N)$ be defined as
\begin{align}
    r_{\Omega}(N) := \sum_{n_1+n_2+n_3=N} \Omega(n_1)\Omega(n_2)\Omega(n_3).
\end{align}
This corresponds to the number of ways of writing an integer $N$ as the sum of three integers $n_1, n_2$ and $n_3$ with solutions weighted based on the number of the total prime divisors of $n_1, n_2$ and $n_3$. We may imagine that the solution is drawn from a multiset wherein each number $n$ occurs with multiplicity $\Omega(n)$.

\begin{theorem} \label{thm:partition theorem 1}
    For all sufficiently large $N$ and every $A>0$ and $\Delta\in(0,\frac{1}{2})$,
\begin{align}
        r_\Omega(N) = \mathfrak{S}\left(N,\left\lfloor A + 4\Delta^{-1}(A+4)\right\rfloor\right)\frac{N^2}{2} + O\bigg(\frac{N^2 (\log \log N)^2}{(\log N)^A}\bigg),
    \end{align}
    where the implied constant depends only on $A$ and $\Delta$. Here
    \begin{align} \label{eq: theorem singular series}
    \mathfrak{S}(N,M) &= \sum_{q = 1}^\infty  \bigg(\frac{\mathfrak{f}(q;N,M)}{N}\bigg)^3 c_q(N),
    \end{align}
    is the singular series associated to $r_\Omega$ and 
    \begin{align} \label{eq:definition of mathfrak f}
         \mathfrak{f}(q;N,M) =
         \frac{N}{q}\bigg((\Omega * \mu)(q) + \sum_{g|q}\frac{g\mu(g)}{\varphi(g)}\sum_{j=1}^M \frac{P_{j,g}(\log \log \frac{Ng}{q})}{(\log \frac{Ng}{q})^{j-1}}\bigg).
    \end{align}
    Here $c_q(n)$ is the Ramanujan sum and $P_j$'s are computable polynomials of degree at most $1$.
\end{theorem}

As a consequence of Theorem~\ref{thm:partition theorem 1} and the Selberg-Delange analysis carried out in Section~\ref{sec:SD_coefficients}, we obtain the following explicit asymptotic.

\begin{corollary} \label{cor:rOmega_asymptotic}
For all sufficiently large $N$,
\begin{align}
    r_\Omega(N) = \frac{N^2}{2}(\log\log N + B(1,1))^3 + O(N^2),
\end{align}
where $B(1,1) = \gamma + \sum_p\bigl[\log(1-1/p) + 1/(p-1)\bigr] = 1.0346534\ldots$ and $\gamma$ is the Euler--Mascheroni constant. In particular, $r_\Omega(N) \sim \frac{N^2}{2}(\log\log N)^3$ and $r_\Omega(N) > 0$ for all sufficiently large $N$.
\end{corollary}

\subsection{Further exponential sum estimates involving $\omega$}
We conclude with two additional results that go beyond the class $\mathcal{F}_0$.
Theorem~\ref{thm:omega squared exp sum} concerns the weight $\omega(n)^2$,
which is no longer additive, while Theorem~\ref{thm:omega-quadratic}
retains the additive weight $\omega(n)$ but replaces the linear phase
$\e(\alpha n)$ with the quadratic phase $\e(\alpha n^2)$.

\begin{theorem} \label{thm:omega squared exp sum}
    Let $\alpha \in \R$, $a \in \Z$, $q \in \N$ and $\Upsilon > 0$ such that
    $|\alpha - \frac{a}{q}| \le \frac{\Upsilon}{q^2}$ with $(a,q) = 1$.
    For any $X \ge 2$, one has
    \begin{align}
        \sum_{n \le X} \omega(n)^2 \e(\alpha n)
        \ll \bigg( \frac{X}{q^{1/5}} \max\{1, \Upsilon^{1/5}\}
        + X^{7/8} + X^{4/5} q^{1/5} \bigg) (\log X)^4.
    \end{align}
\end{theorem}

\begin{theorem}\label{thm:omega-quadratic}
Let $\alpha \in \R$, $a \in \Z$ and $q \in \N$ such that
$|\alpha - \frac{a}{q}| \le q^{-2}$ with $(a,q) = 1$.
For any $X \ge 3$, one has
\begin{equation}\label{eq:main}
\sum_{n\le X}\omega(n)\e(\alpha n^2)
\ll
\frac{X\log\log X}{\sqrt{q}}
+ X\sqrt{\log\log X}
+ (\log\log X)\bigl(\sqrt{X\log q}+\sqrt{q\log q}\,\bigr),
\end{equation}
where the implied constant is absolute.
\end{theorem}

\begin{remark}
The trivial bound is $\sum_{n \le X}\omega(n) \ll X\log\log X$. For $q \ge \log\log X$, the first term in \eqref{eq:main} is $O(X\sqrt{\log\log X})$, and for $1 \le q \le X$ the last two terms are $o(X\log\log X)$, so the bound simplifies to
\[
\sum_{n\le X}\omega(n)\e(\alpha n^2) \ll X\sqrt{\log\log X}
\qquad (q \ge \log\log X,\; q \le X).
\]
This saves a factor of $\sqrt{\log\log X}$ over the trivial estimate.
\end{remark}

\subsection{Structure of the paper} \label{sec:structure}
In Section \ref{sec:exp_sums_additive} we prove Theorem \ref{thm:exp sum bound class F} using elementary tools that gravitate around \eqref{eq: Vinogradov} and bypass techniques from probabilistic number theory and sieve theory. Next, in Section \ref{sec:partitions_1} we prove Theorem \ref{thm:partition theorem 1}. This will require some preliminary lemmas that we will establish along the way such as the convergence of the singular series $\mathfrak{S}(N,M)$ in \eqref{eq: theorem singular series} as well as the Siegel-Walfisz analogue of the summatory function $\sum_{n \le X} \Omega(n)$ over arithmetic progressions. This last step will require revisiting the work of Selberg \cite{SelbergOmega} and Delange \cite{Delange}. Section \ref{sec: proof last 2 theorems} contains the proofs of Theorem \ref{thm:omega squared exp sum} and of Theorem \ref{thm:omega-quadratic}. We conclude with some additional ideas for future work in Section \ref{sec:conclusion}, including an analogous result to Theorem \ref{thm:exp sum bound class F} for a larger class of additive functions and a variety of partition problems that could be tackled with these exponential sums.

\subsection{Notation} 
Throughout the paper, the expressions $f(X)=O(g(X))$, $f(X) \ll g(X)$, and $g(X) \gg f(X)$ are equivalent to the statement that $|f(X)| \le (\ge) C|g(X)|$ for all sufficiently large $X$, where $C>0$ is an absolute constant. A subscript of the form $\ll_{\alpha}$ means the implied constant may depend on the parameter $\alpha$.
The Latin character $p$ will always denote a prime and $\mathbb{P}$ the set of primes.
The prime zeta function is $P(s)=\sum_{p \in \mathbb{P}}p^{-s}$ for $\real(s)>1$. We use the common shorthand $e(\theta)=e^{2\pi i \theta}$, and
$c_q(n)= \sum_{\substack{1\le a\le q \\ (a,q)=1}} e(\frac{a n}{q})$ denotes the Ramanujan sum.
Finally, we shall use the convention that $\delta$, $\varepsilon$, and $\eta$ denote arbitrarily small positive quantities that may not be the same at each occurrence.

\section{Proof of the main result} \label{sec:exp_sums_additive}
%%%%%%%%%%%%%%%%%%%%%%%%%%%%%%%%%%%%%%%%%%%%%%%%%%%%%%%%
In this section, we prove Theorem \ref{thm:exp sum bound class F}.  We start with an auxiliary lemma that could be useful outside the scope of the paper.
\begin{lemma} \label{lem:general principle}
Let $G$ be a non-negative real function, and suppose that there exists $\gamma>0$, $x\ge 1$, $y\ge 1$, such that $G$ satisfies
\begin{equation}\label{gp assumption}
G(\alpha)\ll xq^{-\gamma} + y + x^{1-\gamma}q^{\gamma}
\end{equation}
 whenever $|\alpha-a/q|\le q^{-2}$ and $(q,a)=1$ with $q\in\mathbb N$ and $a\in\mathbb Z$.
Then $G$ also satisfies
\begin{equation}\label{gp conclusion}
G(\alpha)\ll  x(q+x|\alpha q-a|)^{-\gamma}+y+ x^{1-\gamma}(q+x|\alpha q-a|)^{\gamma}
\end{equation}
for all such $\alpha, a, q$.  The implicit constants in \eqref{gp assumption} and \eqref{gp conclusion} are uniform with respect to $x, y, \gamma, \alpha, a, q$.
\end{lemma}

\begin{proof}
Let $\alpha \in \mathbb R$ and fix $a \in \mathbb Z$ and $q\in \mathbb N$ such that  $(q,a)=1$ and $|\alpha-a/q|\le q^{-2}$.  
First observe that if $\alpha=a/q$, then the result is immediate.  Thus we can suppose that $\alpha\not=a/q$.  Now choose $b \in \mathbb Z$ and $ r\in \mathbb N$ such that  $(r,b)=1$, and 
\[
\left|
\alpha-\frac{b}{r}
\right|\le \frac{|\alpha q-a|}{2r} \quad \textnormal{and} \quad r\le \frac{2}{|\alpha q-a|}.
\]
We cannot have $a/q=b/r$, for then we would have $|\alpha-a/q|=|\alpha -b/r|=|\alpha-a/q|/2$ and so $\alpha=a/q$ which is expressly excluded.  Thus
\[
\frac{1}{qr}\le \left|
\frac{b}{r}-\frac{a}{q}
\right|\le \left|
\alpha-\frac{a}{q}
\right| + \left|
\alpha-\frac{b}{r}
\right| \le \left|
\alpha-\frac{a}{q}
\right| + \frac{q}{2r} \left|
\alpha-\frac{a}{q}
\right|\le \left|
\alpha-\frac{a}{q}
\right| + \frac{1}{2rq}.
\]
Hence
\[
\frac{1}{2|\alpha q-a|} \le r\le \frac{2}{|\alpha q-a|}.
\]
Thus by \eqref{gp assumption} we have
\begin{align*}
G(\alpha) &\ll \min\left\{
xq^{-\gamma}+y +x^{1-\gamma}q^{\gamma}, xr^{-\gamma}+y +x^{1-\gamma}r^{\gamma}
\right\}\\
 &\ll \min\left\{
xq^{-\gamma}+y +x^{1-\gamma}q^{\gamma}, x^{1-\gamma}(x|\alpha q-a|)^{\gamma} + y + x(x|\alpha q-a|)^{-\gamma}
\right\}\\
&\ll \min\left\{
xq^{-\gamma}, x(x|\alpha q-a|)^{-\gamma}
\right\} + y + x^{1-\gamma}q^{\gamma} + x^{1-\gamma}(x|\alpha q-a|)^{\gamma}\\
&\ll x(q+x|\alpha q-a|)^{-\gamma} + y + x^{1-\gamma}(q+x|\alpha q-a|)^{\gamma},
\end{align*}
which is what we wanted to prove.
\end{proof}

Equipped with this preliminary result, we may now prove our main result on exponential sums.

\begin{proof}[Proof of Theorem \textnormal{\ref{thm:exp sum bound class F}}]We have 
\begin{equation}
S_f(\alpha; X) = \sum_{n\le X}f(n)\e(\alpha n) = \sum_{\substack{p^\ell m \le X \\ p\, \nmid \, m}}f(p^{\ell}) e(\alpha p^\ell m).
\end{equation}
Let $M\le X$ be a parameter that will be defined later. First consider the terms with $m>X/M$, so that $p^\ell \le M$.  Then by Lemma 2.2 of \cite{Vaughan1997}, together with Lemma \ref{lem:general principle}, we have   
\begin{align}
\bigg|\sum_{p^\ell \le M}  \sum_{\substack{\frac{X}{M} < m \le \frac{X}{p^\ell} \\ p\, \nmid \, m}} f(p^{\ell}) \e(\alpha p^\ell m) \bigg|
&\ll F_f(M)\sum_{p^\ell \le M} \bigg|  \sum_{\substack{\frac{X}{M} < m \le \frac{X}{p^\ell} \\ p\, \nmid \, m}}  \e(\alpha p^\ell m) \bigg| \nonumber \\ 
 & \ll F_f(M)\sum_{p^\ell \le M} \bigg( \bigg|  \sum_{\frac{X}{M} < m \le \frac{X}{p^\ell} }  \e(\alpha p^\ell m) \bigg| + \bigg|  \sum_{\frac{X}{Mp} < k \le \frac{X}{p^{\ell+1}} }  \e(\alpha p^{\ell+1} k) \bigg|\bigg) \label{p divides m} \\ 
 & \ll  F_f(M) \bigg( \sum_{p^\ell \le M}\min\left(\frac{X}{p^{\ell}},\|\alpha p^{\ell}\|^{-1}\right) + \sum_{p^{\ell+1} \le M}\min\left(\frac{X}{p^{\ell+1}},\frac{1}{\|\alpha p^{\ell+1}\|}\right) \nonumber \\ 
 & \quad + \sum_{p \le M^{2/3}}  \bigg|  \sum_{\substack{\frac{X}{Mp} < k \le \frac{X}{p^{\ell+1}} \\ \ell=\lfloor\log_p(M)\rfloor }}  \e(\alpha p^{\ell+1} k) \bigg| 
 + \sum_{M^{2/3} < p \le  M}  \bigg|  \sum_{\frac{X}{Mp} < k \le \frac{X}{p^{2}}}  \e(\alpha p^{2} k) \bigg|\bigg) \label{split off large p}\\  
 & \ll  F_f(M) \bigg( \sum_{n\le M}\min(Xn^{-1},\|\alpha n\|^{-1}) +  \frac{X}{M}\pi( M^{\frac{2}{3}}) +  \frac{X}{M^{\frac{4}3}}\pi\left(M\right) \bigg) \nonumber \\
&\ll F_f(M) \bigg[\left( \frac{X}{q+X|\alpha q-a|}+M+q+X|\alpha q-a| \right)\log 2X  + \frac{X} 
 {M^{\frac{1}{3}}} \bigg]. \label{large m small p}
\end{align}
The second term in line \eqref{p divides m} accounts for the terms with $p\mid m$. The terms in line \eqref{split off large p} represent the contribution from numbers $n=p^{\ell} m$ such that $p|m$ and $p^\ell\le M< p^{\ell+1}$.  In that case, we must have that $\ell=\lfloor\log_p(M)\rfloor$.  In the final term in line \eqref{split off large p}, we necessarily have $\ell=1$ and $p^2 > M^{\frac{4}{3}}$.

Thus it remains to consider
\begin{align}
\sum_{m\le \frac{X}{M}}\sum_{\substack{p^\ell \le \frac{X}{m} \\ p\, \nmid \, m}} f(p^{\ell})  \e(\alpha p^\ell m)
& = \sum_{m\le \frac{X}{M}}\bigg(\sum_{\substack{p \le \frac{X}{m} \\ p\, \nmid \, m}}  \e(\alpha p m) +\sum_{\substack{p^\ell \le \frac{X}{m} \\ p\, \nmid \, m \\ \ell\ge 2}} f(p^{\ell})  \e(\alpha p^\ell m)\bigg) \nonumber \\
& =  \sum_{m\le \frac{X}{M}}( L_1(m) + L_2(m)) .
\end{align}
Notice that
\begin{align}
 \sum_{m\le \frac{X}{M}} L_2(m) \ll  \sum_{m\le \frac{X}{M}} \sum_{\substack{p^\ell \le \frac{X}{m} \\ p\, \nmid \, m \\ \ell\ge 2}} f(p^{\ell})  \ll F_f(X)  \sum_{m\le \frac{X}{M}} \sum_{\substack{p^\ell \le \frac{X}{m} \\ \ell\ge 2}} 1     
\ll XM^{-\frac{1}{2}}F_f(X) (\log X)^{-1}.\label{L_2 est}
\end{align}

We now turn our attention to $L_1$.  We have,
\begin{equation} \label{L_1 split}
\sum_{\substack{p \le \frac{X}{m} \\ p\, \nmid \, m}}  \e(\alpha p m)
 =  \sum_{p \le \frac{X}{m}} \e(\alpha p m) + O\bigg(\bigg| \sum_{\substack{p \le \frac{X}{m} \\ p\, \mid \, m}}  \e(\alpha p m)\bigg|\bigg).
\end{equation}

 For the big-$O$ term of \eqref{L_1 split}, notice that if $p\mid m$ we can write $m = ph$, and the restriction $p\le X/m$ becomes $p^2 \le X/h$.  Summing over $m$, we thus obtain
\begin{equation}\label{L_1 divisors of m}
 \sum_{m\le \frac{X}{M}}\bigg|\sum_{\substack{p \le \frac{X}{m} \\ p\, \mid \, m}} \e(\alpha p m) \bigg| \le  \sum_{h\le \frac{X}{M}}\bigg| \sum_{p^2 \le \frac{X}{h} }  \e(\alpha p^2 h) \bigg|\ll \sum_{h\le \frac{X}{M}} \left(\frac{X}{h}\right)^{1/2}\left(\log\frac{X}{h}\right)^{-1}. 
\end{equation}

Putting together \eqref{L_2 est}, \eqref{L_1 split}, and \eqref{L_1 divisors of m}, we see that 
\begin{equation}\label{reduction to primes}
\sum_{m\le \frac{X}{M}}\sum_{\substack{p^\ell \le \frac{X}{m} \\ p\, \nmid \, m}} f(p^\ell) \e(\alpha p^\ell m) 
= \sum_{m\le \frac{X}{M}}\sum_{p \le \frac{X}{m}} \e(\alpha p m)  +  O\big(XM^{-\frac{1}{2}}F_f(X) (\log X)^{-1}\big).
    \end{equation}

To estimate the main term of \eqref{reduction to primes}, we fix $m\le X/M$ and consider the inner sum.  Let $R =  (\frac{X}{m})^\theta$ for a convenient $\theta \in (0,1)$. Choose $r\in \mathbb N$ and $b \in \mathbb Z$ such that  $ r\le R$, $(b,r)= 1$, and $|\alpha mr -b | \le R^{-1}$.  Then, by Theorem 3.1 of \cite{Vaughan1997}, partial summation, and Lemma \ref{lem:general principle}, we have
\begin{align*} 
  \sum_{p \le \frac{X}{m}} \e(\alpha m p) &\ll (\log X)^3 \bigg(\frac{Xm^{-1}}{(r+\frac{X}{m}|\alpha m r-b|)^{1/2}} +\left(\frac{X}{m}\right)^{\frac{4}{5}} +\left(\frac{X}{m}\right)^{\frac{1}{2}}\left(r+\frac{X}{m}|\alpha m r-b|\right)^{\frac{1}{2}} \bigg) \\
& \ll (\log X)^3\bigg( \frac{Xm^{-1}}{(r+\frac{X}{m}|\alpha m r-b|)^{1/2}}+ \left(\frac{X}{m}\right)^{\frac{4}{5}} + \left(\frac{X}{m}\right)^{\frac{\theta}{2}+ \frac{1}{2}} +  \left(\frac{X}{m}\right)^{1 - \frac{\theta}{2}}\bigg).
\end{align*}
If $r+\frac{X}{m}|\alpha m r-b| \ge \left(\frac{X}{m}\right)^{\frac{2}{5}}$, we have that
$$\sum_{p \le \frac{X}{m}}  \e(\alpha m p) \ll (\log X)^3\bigg(  \bigg(\frac{X}{m}\bigg)^{\frac{4}{5}} + \bigg(\frac{X}{m}\bigg)^{\frac{\theta}{2}+ \frac{1}{2}} +  \bigg(\frac{X}{m}\bigg)^{1 - \frac{\theta}{2}}\bigg).$$
This is optimized when $\frac{2}{5} \le \theta \le \frac{3}{5}$, in which case the right hand side simplifies to $(\log X)^3(\frac{X}{m})^{\frac{4}{5}}$.  The total contribution to \eqref{reduction to primes} in this case is at most
\begin{equation} 
\label{large denominator contribution}
\sum_{m \le \frac{X}{M}} (\log X)^3 \left(\frac{X}{m}\right)^{\frac{4}{5}} \ll  (\log X)^3 X M^{-\frac{1}{5}} . 
\end{equation}

It remains to consider the case when $r+\frac{X}{m}|\alpha m r-b| < \left(\frac{X}{m}\right)^{\frac{2}{5}}$. Fix $a', q'$ such that $1\le q'\le X^{\frac{1}{2}}$, $(a',q')=1$, and $|\alpha-\frac{a'}{q'}| \le (q')^{-1}X^{-\frac{1}{2}}$. If $\frac{b}{mr} \ne \frac{a'}{q'}$, then 
$$\frac{1}{q'mr}\le \left| \frac{b}{mr}-\frac{a'}{q'} \right|\le \left| \alpha-\frac{b}{mr} \right|+\left| \alpha-\frac{a'}{q'} \right|
 \le \frac{1}{mr}\left(\frac{X}m\right)^{-\frac{3}{5}}+\frac{1}{q' X^{\frac{1}2}}. $$
 If we choose $M=4X^{\frac{5}{6}}$, we arrive at the contradiction:
\[1\le  q'\left(\frac{X}m\right)^{-\frac{3}{5}}+\frac{mr}{ X^{\frac{1}2}}\le 2X^{\frac{1}{2}}\left(\frac{X}m\right)^{-\frac{3}{5}}\le 2X^{\frac{1}{2}}M^{-\frac{3}{5}}  = \frac{2}{4^{3/5}}  <  1.\]  
Thus we must have that $bq' = a'mr$.  Since $(r,b)= (q',a') = 1$, we have $r\mid q'$ and $(\frac{q'}{r}, a') = 1$, which gives $\frac{q'}{r}\mid m$.  Write $m = c\cdot\frac{q'}{r}$.  Then 
$$r+\frac{X}{m}|\alpha m r-b| = \frac{r}{q'}\left(q'+ \frac{Xq'}{mr}|\alpha cq' - ca'|\right) = \frac{r}{q'}\left(q'+ X|\alpha q' - a'|\right).$$
Thus the total contribution from terms with $r+\frac{X}{m}|\alpha m r-b| < \left(\frac{X}{m}\right)^{\frac{2}{5}}$ is
\begin{align}
& \ll \sum _{r\mid q'} \sum_{c \le \frac{Xr}{Mq'}} (\log X)^3 \bigg(\left(\frac{Xr}{cq'}\right)\left( \frac{r}{q'}\left(q'+ X|\alpha q' - a'|\right)\right)^{-\frac{1}{2}} +\left(\frac{Xr}{cq'}\right)^{\frac{4}{5}}\bigg) \nonumber \\
 & = \sum_{r\mid q'} \sum_{c \le \frac{Xr}{Mq'}} (\log X)^3 \bigg(\bigg(\frac{X}{c}\bigg)\bigg(\frac{r}{q'}\bigg)^{\frac{1}{2}}\left(q'+X|\alpha q'-a'|\right)^{-\frac{1}{2}} +\left(\frac{Xr}{cq'}\right)^{\frac{4}{5}}\bigg) \nonumber \\ 
& \ll \bigg( \frac{ X\log X}{\left(q'+X|\alpha q'-a'|\right)^{\frac{1}2}} \sum_{r\mid q'} \left(\frac{r}{q'}\right)^{\frac{1}{2}} + \sum_{r\mid q'}\left(\frac{Xr}{q'}\right)^{\frac{4}{5}}\left(\frac{Xr}{Mq'}\right)^{\frac{1}{5}}\bigg) (\log X)^3 \nonumber \\
& = \bigg( \frac{ X(\log X) \sigma_{-1/2}(q')}{\left(q'+X|\alpha q'-a'|\right)^{\frac{1}2}} +  X M^{-\frac{1}{5}}\sigma_{-1}(q')\bigg) (\log X)^3  \nonumber \\ 
& \ll \bigg(\frac{ X}{\left(q'+X|\alpha q'-a'|\right)^{\Delta}} +  X M^{-\frac{1}{5}}\bigg) (\log X)^4 , \label{small denominator contribution}
\end{align}
for any $\Delta\in(0,\frac{1}{2})$.  Here we used that $ \sum_{r\mid q'} (\frac{r}{q'})^{\frac{1}{2}} = \sigma_{-1/2}(q') \ll X^\varepsilon$ and $ \sum_{r\mid q'}\frac{r}{q'} = \sigma_{-1}(q')  \ll \log X$.

Putting \eqref{large denominator contribution} and \eqref{small denominator contribution} into \eqref{reduction to primes} and combining that with \eqref{large m small p}, we have 
\begin{align} S_f(\alpha;X) &  \ll(\log X)^4 \bigg(\frac{ X}{\left(q'+X|\alpha q'-a'|\right)^{\Delta}}  +  X M^{-\frac{1}{5}}\bigg)  + XM^{-\frac{1}{2}}F_f(X) (\log X)^{-1} \nonumber \\
& \quad + F_f(M) \left( \frac{X}{q+X|\alpha q-a|}+M+q+X|\alpha q-a| \right)\log(2X) + F_f(M)XM^{-\frac{1}{3}}\nonumber \\
& \ll\left((\log X)^4 + (\log X)F_f(X)\right)\bigg(\frac{ X}{\left(q'+X|\alpha q'-a'|\right)^{\Delta}} +  X^{\frac{5}{6}}\bigg), \label{intermediate result}
\end{align}
again choosing $M= 4X^{\frac{5}{6}}$.

Now let $a,q$ be as in the statement of the theorem.  If $q\le 2q'$, then we have immediately
\begin{equation} S_f(\alpha;X) \ll \left((\log X)^4 + (\log X)F_f(X)\right)(X(q)^{-\Delta} +X^{\frac{5}{6}}) \label{small q}.\end{equation}
If $q>2q'$, then $\frac{a}{q} \ne \frac{a'}{q'}$ and we have
$$ \frac{1}{q'q} \le \left| \frac{a}{q} - \frac{a'}{q'} \right| \le  \left| \alpha - \frac{a}{q} \right| + \left| \alpha - \frac{a'}{q'} \right| < \frac{1}{q^2} -  \left| \alpha - \frac{a'}{q'} \right| < \frac{1}{2q'q} +  \left| \alpha - \frac{a'}{q'} \right|.$$
Hence 
$ \left| \alpha q' - a' \right| >  \frac{1}{2q}$,
which gives 
\begin{align} \label{large q}
S_f(\alpha;X) & \ll \left((\log X)^4 + (\log X)F_f(X)\right)\bigg(\frac{ X}{\left(Xq^{-1}\right)^{\Delta}} +  X^{\frac{5}{6}}\bigg) \nonumber \\ 
& = \left((\log X)^4 + (\log X)F_f(X)\right)\big(  X^{\frac{5}{6}} + X^{1-\Delta}q^{\Delta}\big). 
\end{align}
Combining \eqref{small q} and \eqref{large q} gives the result.
\end{proof}

\begin{remark}
    The latest improvement on Vinogradov's bound shows that the power of $\log X$ in \eqref{eq: Vinogradov} can be reduced from $3$ to $\frac{5}{2}$ and therefore one could slightly improve the bound for $S_f$ by removing a small fraction from $\log X$. However, as discussed in the introduction, this usually does not matter much as saving powers of $X$ is more valuable than saving powers of $\log X$.
\end{remark}

\begin{remark} \label{remark with Upsilon}
    Another generalization can be accomplished following the technique put forward in \cite{DongRoblesZaharescuZeindler, DongRoblesZaharescuZeindler2} and it would be as follows. Let $\alpha \in \R$, $a \in \Z$ and $q \in \N$ such that $(a,q)=1$ and $|\alpha - a/q| \le \frac{\Upsilon}{q^2}$ with $\Upsilon >0$. Then under the conditions of Theorem \ref{thm:exp sum bound class F} one has
    \begin{align}
        S_f(\alpha; X) \ll \bigg(\frac{X \max\{1,\Upsilon^\Delta\}}{q^\Delta} + X^{\frac{5}{6}}+X^{1-\Delta}q^\Delta\bigg)\Big((\log X)^4 + (\log X)F_f(X)\Big)
    \end{align}
    for any $\Delta \in (0,\frac{1}{2})$. Theorem \ref{thm:exp sum bound class F} follows by taking $\Upsilon=1$.
\end{remark}

The justification for the remark follows by the Dirichlet approximation argument as in \cite[Lemma 2.1]{DongRoblesZaharescuZeindler}. If $\Upsilon \le 1$, the bound follows directly from Theorem \ref{thm:exp sum bound class F} since $\max\{1, \Upsilon^\Delta\} = 1$. If $\Upsilon > 1$, Dirichlet's approximation theorem provides $a_1 \in \Z$ and $q_1 \in \N$ with $(a_1, q_1) = 1$, $q_1 \le 2q$, and $|\alpha - a_1/q_1| \le 1/q_1^2$. If $a_1/q_1 = a/q$, we apply Theorem \ref{thm:exp sum bound class F} directly. Otherwise, the triangle inequality gives $q_1 \ge q/(2\Upsilon)$, and applying Theorem \ref{thm:exp sum bound class F} with denominator $q_1$ yields the stated bound with $\Upsilon^\Delta$ appearing only in the first term.

%%%%%%%%%%%%%%%%%%%%%%%%%%%%%%%%%%%%%%%%%%%%%%%%%%%%%%%%
\section{The ternary Goldbach-Vinogradov problem associated to Omega} \label{sec:partitions_1}
In this section we shall use the circle method to prove Theorem \ref{thm:partition theorem 1}, following a blueprint similar to that of \cite[$\mathsection$ 8]{Nathanson}.  Throughout this section, we adopt the shorthand notation
\begin{align}
    F_N(\alpha) := \sum_{n \le N} \Omega(n) \e(\alpha n) \quad \textnormal{for} \quad S_\Omega(\alpha;N).
\end{align}

\subsection{Preliminary lemmas}
A key ingredient is the Selberg-Delange method, which provides asymptotic expansions for sums of the form $\sum_{n \le x} z^{\Omega(n)}$ (and variants over arithmetic progressions) in powers of $(\log x)^{z-1}$, uniformly for $z$ in a neighborhood of $1$. The method, originating in the work of Selberg~\cite{SelbergOmega} and Delange~\cite{Delange}, produces coefficients $\lambda_j(z,q)$ that are analytic in $z$; differentiating at $z=1$ recovers asymptotics for $\sum \Omega(n)$. See Tenenbaum~\cite{tenenbaum} for a modern treatment. The following lemma summarizes the expansion we need.
%%%%%%%%%%%
\begin{lemma} \label{lem:SiegelWalfiszOmega}
    Let $\rho$ be defined by
     \begin{align} \label{eq:rho_def}
    \rho(s;z,q) := \prod_{p|q} \bigg(1-\frac{z}{p^s}\bigg).
    \end{align}
    There exists a family $\{z \mapsto A_r(z):r=1,2,3,\cdots\}$ of holomorphic functions inside the open disk $\{z \in \C: |z|<2\}$ such that, for all $j \ge 1$, with
     \begin{align}
    \lambda_j(z,q) := \sum_{r=1}^j \binom{z-r}{j-r} A_r(z) \rho^{(j-r)}(1;z,q),
    \end{align}
    one has for $|z| \le 2-\delta$, $x \ge 1$, $1 \le q \le (\log x)^b$, for some $b>0$, $(h,q)=1$ and $t \ge 1$
     \begin{align} \label{eq: siegel Walfizes z omega}
    \sum_{\substack{n < x \\ n \equiv h \modu q}} z^{\Omega(n)}= \frac{x}{\varphi(q)} \sum_{j=1}^t \lambda_j(z,q)(\log x)^{z-j} + O_{\delta,t,b}\bigg(\frac{x}{\varphi(q)} (\log x)^{z-t-1} (\log \log q)^{t+2}\bigg).
    \end{align}
    Moreover, one has the upper bound
     \begin{align}
    |\lambda_j(z,q)| \ll _{\delta,j} (\log \log q)^{j+1}.
    \end{align}
\end{lemma}

%%%%%%%%%%%
\begin{proof}
    See \cite[Theorem 3]{DupainHallTenenbaum1982}.
\end{proof}
%%%%%%%%%
Differentiating with respect to $z$ and then letting $z \to 1$ leads to
 \begin{align}
\sum_{\substack{n < x \\ n \equiv h \modu q}} \Omega(n) = \frac{x}{\varphi(q)} \sum_{j=1}^t \frac{P_{j,q}(\log \log x)}{(\log x)^{j-1}} + O_{t,b} \bigg(\frac{x \log \log x}{\varphi(q) (\log x)^t} (\log \log q)^{t+2}\bigg),
\end{align}
where $P_{j,q}(y)$ are polynomials in $y$ of degree at most one, given by
 \begin{align} \label{eq:Pjg_def}
P_{j,q}(\log \log x): = b(j,q)\log \log x + B(j,q)
\end{align}
with $b$ and $B$ given by
 \begin{align}
b(j,q): = {\lambda _j}(1,q) \quad \textnormal{and} \quad B(j,q): = \mathop {\lim }\limits_{z \to 1} \frac{d}{{dz}}{\lambda _j}(z,q).
\end{align}
The explicit form of the holomorphic functions $A_r$, and hence the explicit form of $\lambda_j$, can be found using the techniques in \cite[$\mathsection$ 14]{IvicBook} as well as \cite[Chapter II.6]{tenenbaum}.
\begin{lemma} \label{lem:FxA/Q}
    Let $\mathfrak{f}$ be as in \eqref{eq:definition of mathfrak f}, let $M$ be an arbitrary positive integer, and  let $N$ be a large positive real number. For $1\le x \le N$, $1 \le q \le (\log x)^{b}$ for some $b>0$, and $(a,q)=1$, we have 
    \begin{align}
         F_x\bigg(\frac{a}{q}\bigg) = \mathfrak{f}(q;x,M) + O_{M} \bigg( \frac{N q^{\varepsilon}\log\log N }{(\log N)^{M}}\bigg),
    \end{align}
     for any $\varepsilon>0$, where $Q = (\log N)^b$.
\end{lemma}

%%%%%%%%%%%
\begin{proof}
We start by writing
\begin{align}
    F_x\bigg(\frac{a}{q}\bigg) = \sum_{n \le x} \Omega(n) \e \bigg(n\frac{a}{q}\bigg) = \sum_{1 \le r \le q} \sum_{\substack{n \le x \\ n \equiv r \modu q}} \Omega(n) \e \bigg(n\frac{a}{q}\bigg) = \sum_{1 \le r \le q} \e \bigg(r\frac{a}{q}\bigg) \sum_{\substack{n \le x \\ n \equiv r \modu q}} \Omega(n)  
\end{align}
where $(a,q)=1$. The inner sum presents a problem since $r$ and $q$ are not necessarily coprime. Let $g=(r,q)$. We may then write the inner sum as
\begin{align}
    \sum_{\substack{n \le x \\ n \equiv r \modu q}} \Omega(n) = \sum_{\substack{1 \le m \le \frac{x}{g} \\ m \equiv \frac{r}{g} \modu \frac{q}{g}}} \Omega(gm) = \sum_{\substack{1 \le m \le \frac{x}{g} \\ m \equiv \frac{r}{g} \modu \frac{q}{g}}} [\Omega(g)+\Omega(m)],
\end{align}
by the total additivity of the function $\Omega$ (i.e., $\Omega(mn)=\Omega(m)+\Omega(n)$ for all $m$ and $n$). Therefore 
\begin{align}
    \sum_{\substack{n \le x \\ n \equiv r \modu q}} \Omega(n) = \Omega(g) \bigg(\frac{x}{q}+O(1)\bigg) + \sum_{\substack{1 \le m \le \frac{x}{g} \\ m \equiv \frac{r}{g} \modu \frac{q}{g}}} \Omega(m).
\end{align}
Inserting this into $F_x$ we see that
\begin{align}
    F_x\bigg(\frac{a}{q}\bigg) =  \sum_{g|q} \sum_{\substack{1 \le r \le q \\ (r,q)=g}} \e \bigg(r\frac{a}{q}\bigg) \bigg[ \Omega(g)\bigg(\frac{x}{q}+O(1)\bigg) + \sum_{\substack{1 \le m \le \frac{x}{g} \\ m \equiv \frac{r}{g} \modu \frac{q}{g}}} \Omega(m)\bigg] =: S_1 + S_2.
\end{align}

For the first sum we have
\begin{align}
    S_1 =  \sum_{g|q} \sum_{\substack{1 \le r \le q \\ (r,q)=g}} \e \bigg(r\frac{a}{q}\bigg) \Omega(g)\bigg(\frac{x}{q}+O(1)\bigg) = \frac{x}{q} \sum_{g|q} \Omega(g) \sum_{\substack{1 \le r \le q \\ (r,q)=g}} \e \bigg(r\frac{a}{q}\bigg) + O\bigg(\sum_{1 \le r \le q}\Omega(r)\bigg).
\end{align}
The exponential sum can be evaluated by writing $r=gr'$ and $q=gq'$ so that $(r,q)=g$ if, and only if, $(r',q')=1$. Therefore
\begin{align} \label{eq: exponential r sum}
    \sum_{\substack{1 \le r \le q \\ (r,q)=g}} \e \bigg(r\frac{a}{q}\bigg) = \sum_{\substack{1 \le r' \le q' \\ (r',q')=1}} \e \bigg(r'\frac{a}{q'}\bigg) = c_{q'}(a) = \mu(q') = \mu \bigg(\frac{q}{g}\bigg),
\end{align}
since $(a,q')=1$. This leaves us with 
\begin{align}
    S_1 = \frac{x}{q}\sum_{g|q} \Omega(g) \mu\bigg(\frac{q}{g}\bigg) + O(q \log \log q) = \frac{x}{q}(\Omega*\mu)(q) +O(q \log \log q).
\end{align}

We now consider the second sum,
\begin{align}
    S_2 = \sum_{g|q} \sum_{\substack{1 \le r \le q \\ (r,q)=g}}\e \bigg(r\frac{a}{q}\bigg) \sum_{\substack{1 \le m \le \frac{x}{g} \\ m \equiv \frac{r}{g} \modu \frac{q}{g}}} \Omega(m).
\end{align}
We now use Lemma \ref{lem:SiegelWalfiszOmega} to write 
\begin{align}
 & \sum_{\substack{1 \le m \le \frac{x}{g} \\ m \equiv \frac{r}{g} \modu \frac{q}{g}}} \Omega(m)
  =\frac{\frac{x}{g}}{\varphi(\frac{q}{g})}\sum_{j=1}^M \frac{P_{j,\frac{q}{g}}(\log \log \frac{x}{g})}{(\log \frac{x}{g})^{j-1}} + O \bigg(\frac{x\log\log x (\log\log q)^{M+2}}{(\log  x)^{M}}\bigg).
\end{align}
Since this estimate is independent of $r$, we may move the sum over $r$ to the inside and evaluate it using \eqref{eq: exponential r sum}.  This yields 
\begin{align}
    S_2 &= \sum_{g|q}\frac{\mu(\frac{q}{g})}{\varphi(\frac{q}{g})}\frac{x}{g}\sum_{j=1}^M \frac{P_{j,\frac{q}{g}}(\log \log \frac{x}{g})}{(\log \frac{x}{g})^{j-1}} + O \bigg( \sum_{g|q} \frac{x \log\log x (\log\log q)^{M+2}}{(\log  x)^{M}}\bigg) \nonumber \\
 &= \frac{x}{q}\sum_{g|q}\frac{g\mu(g)}{\varphi(g)}\sum_{j=1}^M \frac{P_{j,g}(\log \log \frac{xg}{q})}{(\log \frac{xg}{q})^{j-1}} +  O \bigg( \frac{N q^{\varepsilon}\log\log N }{(\log N)^{M}}\bigg)  .
\end{align}
Therefore, adding up the contributions from $S_1$ and from $S_2$ we arrive at
\begin{align}
    F_x\bigg(\frac{a}{q}\bigg) 
     = \frac{x}{q}\bigg((\Omega * \mu)(q) + \sum_{g|q}\frac{g\mu(g)}{\varphi(g)}\sum_{j=1}^M \frac{P_{j,g}(\log \log \frac{xg}{q})}{(\log \frac{xg}{q})^{j-1}}\bigg) 
    + O \bigg( \frac{N q^{\varepsilon}\log\log N }{(\log N)^{M}}\bigg).
\end{align}
The main term on the right-hand side is recognized to be $\mathfrak{f}(q;x,M)$.
\end{proof}
%%%%%%%%%%%
\begin{lemma} \label{lem:f_bound}
     For $M,q$ positive integers and $N$ a real number, we have
    \begin{align}
        \mathfrak{f}(q;N,M)  \ll \frac{NM \log \log N}{q^{1-\varepsilon}}
    \end{align}
    for any $\varepsilon>0$.
\end{lemma}

%%%%%%%%%%%
\begin{proof}
    We write
    \begin{align*}
        \mathfrak{f}(q;N,M) = 
        \frac{N}{q}\bigg((\Omega * \mu)(q) + \sum_{g|q}\frac{g\mu(g)}{\varphi(g)}\sum_{j=1}^M \frac{P_{j,g}(\log \log \frac{Ng}{q})}{(\log \frac{Ng}{q})^{j-1}}\bigg).
    \end{align*}
     Noticing that $\Omega * \mu$ and $\mu$ both take values in $\{-1,0,1\}$ and recalling that $P_{j,q}$ are polynomials of degree at most one, we have 
         \begin{align}
        \mathfrak{f}(q;N,M) 
        &\ll \frac{N}{q}\bigg(1 + M \log \log N\sum_{g|q}\frac{g}{\varphi(g)} \bigg).
        \end{align}
    Since $\varphi(a) > a^{1-\delta}$ and $\sum_{a|n} a^\nu \ll n^{\nu+\delta}$ for every $\delta>0$, we have 
    $$ \sum_{g|q}\frac{g}{\varphi(g)}  \ll \sum_{g|q}g^{\delta}\ll q^{2\delta}.$$
    The result now follows by setting $\delta =\varepsilon/2$.
\end{proof}
%%%%%%%%%%%
\begin{lemma}
    Let $c_q(N)$ denote the Ramanujan sum. The singular series $\mathfrak{S}(N,M)$ defined by
    \begin{align}
    \mathfrak{S}(N,M) &= \sum_{q = 1}^\infty  \bigg(\frac{\mathfrak{f}(q;N,M)}{N}\bigg)^3 c_q(N),
    \end{align}
    converges absolutely. Moreover, for any $\eta>0$ one has that
     \begin{align}
    \mathfrak{S}(N,M;Q) = \sum_{1 \le q \le Q} \bigg(\frac{\mathfrak{f}(q;N,M)}{N}\bigg)^3 c_q(N) = \mathfrak{S}(N,M) + O\bigg(\frac{(M\log \log N)^3}{Q^{1-\eta}}\bigg),
    \end{align}
    where the implied constant depends only on $\eta$.
\end{lemma}

\begin{proof}
    Clearly $|c_q(N)|\le q$. Let $0<\varepsilon<\frac{1}{3}$ be arbitrary. By Lemma \ref{lem:f_bound} there exists a constant $C_{\varepsilon}>0$ so that 
    \begin{align}
        \bigg|\frac{\mathfrak{f}(q;N,M)}{N}\bigg|^3 |c_q(N)|
        \le C_{\varepsilon} \frac{(M\log \log N)^3}{q^{2-3\varepsilon }}.
    \label{f cubed abs bound}
    \end{align}
      This shows that $\mathfrak{S}(N,M)$ converges absolutely by the comparison test. Moreover, we have 
    \begin{align}
    |\mathfrak{S}(N,M)  - \mathfrak{S}(N,M;Q)| 
     &= \bigg| \sum_{q>Q}\bigg(\frac{\mathfrak{f}(q;N,M)}{N}\bigg)^3 c_q(N) \bigg|  \nonumber \\
    &\le C_{\varepsilon} (M\log \log N)^3  \sum_{q > Q} \frac{1}{q^{2-3\varepsilon }} \nonumber 
     \ll \frac{(M\log \log N)^3}{Q^{1-3\varepsilon}},
    \end{align}
    from which the result follows by choosing $\varepsilon = \min\{\frac{\eta}{3}, \frac{1}{6}\}$.
\end{proof}
%%%%%%%%%%%
\begin{lemma} \label{lem:F and F cubed}
    Let $B$ be positive real number and $M$ be a positive integer such that $M +1 > B$. If $1 \le q \le Q = (\log N)^B$ and $(a,q)=1$, then for $\alpha = a/q+\beta$ with $|\beta| \le Q/N$, we have
    \begin{align} \label{eq:FN_power1}
        F_N(\alpha) = \frac{\mathfrak{f}(q;N,M)}{N} u(\beta) + O \bigg(\frac{NQ^{1+\varepsilon}}{(\log N)^{M+1}}\bigg)
    \end{align}
    as well as
    \begin{align} \label{eq:FN_power3}
        F_N(\alpha)^3 = \bigg(\frac{\mathfrak{f}(q;N,M)}{N}\bigg)^3u(\beta)^3 + O \bigg(\frac{N^3Q^{1+\varepsilon}(\log \log N)^2}{(\log N)^{M+1}}\bigg),
    \end{align} 
    where the implied constant depends on $B$ and $M$. Here $u(\beta) = \sum_{n \le N}\e(n \beta)$.
\end{lemma}
%%%%%%%%%%%
\begin{proof}
  We have
     \begin{align}
        F_N(\alpha) - \frac{\mathfrak{f}(q;N,M)}{N} u(\beta) 
        &= \sum_{n \le N} \Omega(n) \e(\alpha n) - \frac{\mathfrak{f}(q;N,M)}{N}\sum_{n \le N} \e(\beta n) \nonumber \\
        &= \sum_{n \le N} \Omega(n) \e\bigg(n\frac{a}{q}\bigg)\e(\beta n) - \frac{\mathfrak{f}(q;N,M)}{N} \sum_{n \le N} \e(\beta n) \nonumber \\
        &= \sum_{n \le N} \bigg\{\Omega(n) \e\bigg(n\frac{a}{q}\bigg) -  \frac{\mathfrak{f}(q;N,M)}{N} \bigg\}\e(\beta n).
    \end{align}
    We apply partial summation to obtain
    \begin{align} 
        F_N(\alpha) - \frac{\mathfrak{f}(q;N,M)}{N} u(\beta) &= A(N) \e(N \beta)- 2 \pi i \beta \int_1^N A(x) e(\beta x)dx \nonumber \\
        &\ll |A(N)|  + |\beta|N\max\{A(x)\,:\, 1\le x\le N\}, \label{eq: part sum setup}
    \end{align}
    where
      \begin{align}
        A(x) := \sum_{n \le x} \bigg(\Omega(n) \e\bigg(n\frac{a}{q}\bigg) -  \frac{\mathfrak{f}(q;x,M)}{x} \bigg). 
    \end{align}
    We now need to bound $A(x)$ for  $1 \le x \le N$.  Lemma \ref{lem:FxA/Q} requires $q\le (\log x)^b$ for some $b$.  Choosing $b=2B$, we split the interval $1\le x \le N$ into $1\le x\le \exp(q^{\frac{1}{2B}})$ and $\exp(q^{\frac{1}{2B}}) \le x \le N$.  
    For small $x$, we have 
    \begin{equation} x\le \exp(q^{\frac{1}{2B}}) \le \exp(Q^{\frac{1}{2B}}) \le \exp(\sqrt{\log N}) \ll N^\varepsilon.
    \end{equation}
    Thus by Lemma \ref{lem:f_bound}, we  have for $1\le x\le \exp(q^{\frac{1}{2B}})$
    \begin{align}
        A(x)  \ll \sum_{n \le x} \Omega(n)  +   \mathfrak{f}(q;x,M)  
         \ll x\log\log x \ll N^\varepsilon.
    \end{align}
    Meanwhile, for $\exp(q^{\frac{1}{2B}}) \le x \le N$, we have by Lemma \ref{lem:FxA/Q} and Lemma \ref{lem:f_bound} that 
    \begin{align}
        A(x) & = \sum_{n \le x} \bigg(\Omega(n) \e\bigg(n\frac{a}{q}\bigg) -  \frac{\mathfrak{f}(q;x,M)}{x} \bigg) \nonumber \\
        &=\sum_{n \le x} \Omega(n) \e\bigg(n\frac{a}{q}\bigg) - \mathfrak{f}(q;x,M) + O\bigg(\frac{\mathfrak{f}(q;x,M)}{x}\bigg)  \nonumber \\
        &=F_x\bigg(\frac{a}{q}\bigg) - \mathfrak{f}(q;x,M) + O_M(q^{\varepsilon} \log \log x) \nonumber 
        \ll_M \frac{NQ^{\varepsilon}}{(\log N)^{M+1}} ,
    \end{align}
    since $\log\log x \le \log \log N = \log Q^{1/B} \ll Q^\varepsilon.$
Putting this into \eqref{eq: part sum setup} yields
\begin{align} 
        F_N(\alpha) - \frac{\mathfrak{f}(q;N,M)}{N} u(\beta) 
        &\ll |A(N)|  + |\beta|N\max\{A(x)\,:\, 1\le x\le N\} \nonumber \\
       &\ll \frac{NQ^{\varepsilon}}{(\log N)^{M+1}} + \frac{|\beta|N^2 Q^{\varepsilon}}{(\log N)^{M+1}} \ll \frac{N Q^{1+\varepsilon}}{(\log N)^{M+1}}.
    \end{align}
    since $|\beta| \le Q/N$. This shows the estimate for $F_N(\alpha)$.

    Next, we note that $|u(\beta)| \le N$ and that the main term in \eqref{eq:FN_power1} is $\ll_M  N \log \log N$. Moreover, the error term in \eqref{eq:FN_power1} holds for all $\varepsilon>0$, in particular, for $\varepsilon=(M+1-B)/(2B)$.  We then have
    \[
    \frac{N Q^{1+\varepsilon}}{(\log N)^{M+1}} = \frac{N}{(\log N)^{M+1-B - \varepsilon B}} < N.
    \]
    This implies that
    \[
    F_N(\alpha)^3 = \bigg(\frac{\mathfrak{f}(q;N,M)}{N}u(\beta)\bigg)^3 + O\bigg(N^2 (\log \log N)^2 \frac{NQ^{1+\varepsilon}}{(\log N)^{M+1}}\bigg)
    \]
    and the estimate for $F_N(\alpha)^3$ now follows.
\end{proof}
%%%%%%%%%%%
\subsection{The major and minor arcs}
%%%%%%%%%%%
We decompose the unit interval $[0,1]$ into two disjoint sets: the major and minor arcs, which we define as follows. Let $B>0$ and $Q = (\log N)^B$. For $1 \le q \le Q$ and $0 \le a \le q$ with $(a,q)=1$ the major arc $\mathfrak{M}(q,a)$ is the interval consisting of all real numbers $\alpha \in [0,1]$ such that $|\alpha - \frac{a}{q}| \le \frac{Q}{N}$. The major arcs are pairwise disjoint for large $N$. The set of major arcs is given by the union
\begin{align}
    \mathfrak{M} := \bigcup_{1 \le q \le Q} \bigcup_{\substack{0 \le a \le q \\ (a,q)=1}} \mathfrak{M} (a,q).
\end{align}
The minor arcs will therefore be $\mathfrak{m}=[0,1] \backslash \mathfrak{M}$. By orthogonality we have that 
\begin{align}
    r_\Omega(N) = \int_0^1 F_N(\alpha)^3\e(-\alpha N)d\alpha = \bigg(\int_{\mathfrak{M}}+\int_{\mathfrak{m}}\bigg)F_N(\alpha)^3\e(-\alpha N)d\alpha.
\end{align}

Our task is to now extract the main term for $r_\Omega(N)$ from the integral over $\mathfrak{M}$ and show that the the integral over $\mathfrak{m}$ is of lower order.  We first use the results from the previous section to compute the integral over the major arcs.
\begin{theorem}[The major arcs] \label{thm:major_arcs} 
    Let B be a positive real number and $M$ a positive integer with $M+1 > B$.   Then the integral over the major arcs is
    \begin{align}
        \int_{\mathfrak{M}} F_N(\alpha)^3 \e(-\alpha N) d\alpha = \mathfrak{S}(N,M)\frac{N^2}{2} + O_{B,M}\bigg(\frac{N^2}{(\log N)^C}\bigg),
    \end{align}
    for any number $C < \min\{B,\  M+1-4B\}$. 
\end{theorem} 
%%%%%%%%%%%
\begin{proof}
    The first observation is that the length of the major arcs $\mathfrak{M}(q,a)$ is $Q / N$ if $q=1$ and $2Q/N$ if $q \ge 2$.
    We use Lemma \ref{lem:F and F cubed} to handle the error term as
    \begin{align}
        E_{\mathfrak{M}} :&= \int_{\mathfrak{M}} \bigg[F_N(\alpha)^3 - \bigg(\frac{\mathfrak{f}(q;N,M)}{N}\bigg)^3 u(\beta)^3\bigg] \e(-N \alpha) d\alpha \nonumber \\
        &= \sum_{1 \le q \le Q} \sum_{\substack{0 \le a \le q \\ (a,q)=1}} \int_{\mathfrak{M}(q,a)} \bigg[F_N(\alpha)^3 - \bigg(\frac{\mathfrak{f}(q;N,M)}{N}\bigg)^3 u(\beta)^3\bigg] \e(-N \alpha) d\alpha \nonumber \\
        &\ll \sum_{1 \le q \le Q} \sum_{\substack{0 \le a \le q \\ (a,q)=1}} \int_{\mathfrak{M}(q,a)}\frac{N^3Q^{1+\varepsilon}(\log \log N)^2}{(\log N)^{M+1}} d\alpha \nonumber \\
        &\ll \sum_{1 \le q \le Q} \sum_{\substack{0 \le a \le q \\ (a,q)=1}} \frac{N^2Q^{2+\varepsilon}(\log \log N)^2}{(\log N)^{M+1}}  \nonumber \\
        &\le \frac{N^2Q^{4+\varepsilon}(\log \log N)^2}{(\log N)^{M+1}}  \nonumber \\
        &\le \frac{N^2(\log \log N)^2}{(\log N)^{M+1-4B-\varepsilon B}},
    \end{align}
    for all $\epsilon>0$.
    On the other hand, if $\alpha = a/q+\beta \in \mathfrak{M}(q,a)$, then $|\beta| \le Q/N$ and we see that for the main term 
    \begin{align}
        M_{\mathfrak{M}} :&= \int_{\mathfrak{M}} \bigg(\frac{\mathfrak{f}(q;N,M)}{N}\bigg)^3 u\bigg(\alpha-\frac{a}{q}\bigg)^3 \e(-N \alpha)d\alpha \nonumber \\
        &=\sum_{1 \le q \le Q} \sum_{\substack{0 \le a \le q \\ (a,q)=1}} \bigg(\frac{\mathfrak{f}(q;N,M)}{N}\bigg)^3 \int_{\mathfrak{M}(q,a)} u\bigg(\alpha-\frac{a}{q}\bigg)^3 \e(-N \alpha)d\alpha \nonumber \\
        &=\sum_{1 \le q \le Q} \bigg(\frac{\mathfrak{f}(q;N,M)}{N}\bigg)^3\sum_{\substack{0 \le a \le q \\ (a,q)=1}} \e\bigg( -\frac{Na}{q}\bigg) \int_{a/q-Q/N}^{a/q+Q/N} u\bigg(\alpha-\frac{a}{q}\bigg)^3 \e(-N \alpha)d\alpha \nonumber \\
        &=\sum_{1 \le q \le Q} \bigg(\frac{\mathfrak{f}(q;N,M)}{N}\bigg)^3 c_q(N) \int_{-Q/N}^{Q/N} u(\beta)^3 \e(-N \beta)d\beta \nonumber \\
        &=\mathfrak{S}(N,M;Q) \int_{-Q/N}^{Q/N} u(\beta)^3 \e(-N \beta)d\beta.
    \end{align}
    The integral over $\beta$ was evaluated in \cite[$\mathsection$ 8.4]{Nathanson} and it was shown to be $\frac{N^2}{2} + O(\frac{N^2}{Q^2})$. Consequently, we arrive at
    \begin{align}
    \int_\mathfrak{M} &F_N(\alpha )^3\e( - N\alpha )d\alpha \nonumber \\
    &= \mathfrak{S}(N,M;Q)\int_{ - Q/N}^{Q/N}  u(\beta )^3\e( - N\beta )d\beta  + O\bigg( \frac{N^2(\log \log N)^2}{(\log N)^{M + 1 - 4B - \varepsilon B}} \bigg) \nonumber \\
    &= \bigg( \mathfrak{S}(N,M) + O\bigg( \frac{(\log \log N)^3}{Q^{1 - \eta}} \bigg) \bigg)\left( \frac{N^2}{2} + O\left( \frac{N^2}{Q^2} \right) \right) + O\bigg( \frac{N^2(\log \log N)^2}{(\log N)^{M + 1 - 4B - \varepsilon B }} \bigg) \nonumber \\
    &= \mathfrak{S}(N,M)\frac{N^2}{2} + O\bigg( \frac{N^2(\log \log N)^3}{Q^{1 - \eta}} \bigg) + O\bigg( \frac{N^2(\log \log N)^2}{(\log N)^{M + 1 - 4B - \varepsilon B}} \bigg) \nonumber\\
    &= \mathfrak{S}(N,M)\frac{N^2}{2} + O\bigg( \frac{N^2(\log \log N)^3}{(\log N)^{B(1 - \eta)}} \bigg) + O\bigg( \frac{N^2(\log \log N)^2}{(\log N)^{M + 1 - 4B - \varepsilon B}} \bigg). \label{sing series error proof}
    \end{align}
   Note that $(\log\log N)^3 \ll (\log N)^\delta$ for all $\delta>0$.  Thus, choosing $\varepsilon$ and $\eta$ sufficiently small in \eqref{sing series error proof}, we have 
   \begin{align}
    \int_\mathfrak{M} &F_N(\alpha )^3\e( - N\alpha )d\alpha \nonumber 
      = \mathfrak{S}(N,M)\frac{N^2}{2} + O\bigg( \frac{N^2}{(\log N)^{C}} \bigg),
    \end{align}
    and this ends the proof.
\end{proof}
%%%%%%%%%%%
We can encapsulate the application of Theorem \ref{thm:exp sum bound class F} into the following result for the minor arcs.
%%%%%%%%%%%
\begin{theorem}[The minor arcs] \label{thm:minor_arcs}
    For any $B>0$, we have
    \begin{align}
        \int_{\mathfrak{m}} F_N (\alpha)^3 \e(-\alpha N) d\alpha \ll \frac{N^2 (\log \log N)^2}{(\log N)^{B\Delta-4}}
    \end{align}
    for any $\Delta \in (0,\frac{1}{2})$ and where the implied constant depends only $B$. 
\end{theorem}
%%%%%%%%%%%
\begin{proof}
    Let $\alpha \in \mathfrak{m}$. By Dirichlet's theorem, for any real $\alpha$ there exists a fraction $a/q \in [0,1]$ with $1 \le q \le N/Q$ and $(a,q)=1$ such that
    \begin{align}
        \bigg|\alpha-\frac{a}{q}\bigg| \le \frac{Q}{qN} \le \min \bigg(\frac{Q}{N},\frac{1}{q^2}\bigg).
    \end{align}
    Moreover, $Q < q \le \frac{N}{Q}$, for otherwise $q \le Q$ and then $\alpha \in \mathfrak{M}(q,a) \subseteq \mathfrak{M}$, which would be a contradiction. By Theorem \ref{thm:exp sum bound class F} with $f = \Omega$ we have
    \begin{align}
        F_N(\alpha) \ll \bigg(\frac{N}{q^{\Delta}}+N^{5/6}+N^{1-\Delta}q^{\Delta}\bigg) (\log N)^4.
    \end{align}
    for any $\Delta \in (0, \frac{1}{2})$. The above lower and upper bounds on $q$ yields
        \begin{align}
        F_N(\alpha) \ll \bigg(\frac{N}{(\log N)^{B\Delta}}+N^{5/6}+N^{1-\Delta}\bigg( \frac{N}{\log N}\bigg)^{\Delta}\bigg) (\log N)^4
        \ll \frac{N}{(\log N)^{B\Delta-4}}.
    \end{align}
    Now, we have by \cite[Eq. (14.146)]{IvicBook} that 
    \begin{align}
        \int_0^1 |F_N(\alpha )|^2d\alpha = \sum_{n \le N} \Omega (n)^2 \le N (\log \log N)^2.
    \end{align}
    Consequently, we end up with
    \begin{align}
        \int_{\mathfrak{m}} |F_N(\alpha)|^3 d\alpha \ll \sup\{|F_N(\alpha)| : 
        \alpha \in \mathfrak{m}\} \int_{\mathfrak{m}} |F_N(\alpha)|^2 d\alpha \nonumber &\ll \frac{N}{(\log N)^{B\Delta-4}} \int_0^1 |F_N(\alpha)|^2 d\alpha \nonumber \\
        &\ll \frac{N^2 (\log \log N)^2}{(\log N)^{B\Delta-4}},
    \end{align}
    which was the last ingredient of the proof.
\end{proof}
%%%%%%%%%%%
\subsection{The asymptotic formula}
%%%%%%%%%%%

We now prove Theorem \ref{thm:partition theorem 1} by 
    merging our results from the major arcs (Theorem \ref{thm:major_arcs}) and minor arcs (Theorem \ref{thm:minor_arcs}). For any positive number $B$ and positive integer $M$ with $M+1 > 2B$, and any $\Delta \in (0, \frac{1}{2})$, we have 
    \begin{align}
        r{_\Omega} (N) &=\bigg(\int_{\mathfrak{M}}+\int_{\mathfrak{m}}\bigg) F_N(\alpha)^3 \e(-\alpha N)d\alpha \nonumber \\
        &=\mathfrak{S}(N,M)\frac{N^2}{2} + O\bigg( \frac{N^2}{(\log N)^{C}} \bigg) +  O \bigg( \frac{N^2 (\log \log N)^2}{(\log N)^{B\Delta-4}} \bigg), 
    \end{align}
where $C$ is any number satisfying $C<\min\{B,\ M+1-4B\}$.  

Let $A>0$ and $\Delta\in(0,1)$ be given.  Choose $B=\Delta^{-1}(A+4)$ and $M=\lfloor A + 4B\rfloor$.  Then $M+1 > B$ and 
$$\min\{B,\ M+1-4B,\ B\Delta-4\} = B\Delta-4 = A.$$
Thus we arrive at 
  \begin{align}
        r_\Omega (N) 
        =\mathfrak{S}\left(N,\left\lfloor A + 4\Delta^{-1}(A+4)\right\rfloor\right)\frac{N^2}{2} +  O \bigg(\frac{N^2 (\log \log N)^2}{(\log N)^{A}} \bigg), 
    \end{align}
    as desired.

\subsection{Explicit Selberg-Delange coefficients} \label{sec:SD_coefficients}
We now make the polynomials $P_{j,g}$ from \eqref{eq:Pjg_def} explicit. This will be needed in Section~\ref{sec:positivity} to determine the scale and positivity of the singular series $\mathfrak{S}(N,M)$.

\begin{lemma}[Coefficients for $j=1$] \label{lem:j1_coefficients}
Let $P_{1,g}(y) = b(1,g)\, y + B(1,g)$ be the polynomial from the Selberg-Delange expansion. Then:
\begin{enumerate}
    \item[(i)] $b(1,g) = \dfrac{\varphi(g)}{g}$ for all squarefree $g \geq 1$.
    
    \item[(ii)] $B(1,g) = \dfrac{\varphi(g)}{g}\bigg(B(1,1) - \displaystyle\sum_{p \mid g} \dfrac{1}{p-1}\bigg)$ for all squarefree $g \geq 1$, where
    \begin{align}
        B(1,1) = \gamma + A_1'(1)
    \end{align}
    with $\gamma = 0.5772156649\ldots$ the Euler--Mascheroni constant and
    \begin{align}
        A_1'(1) = \sum_p \left( \log\left(1 - \frac{1}{p}\right) + \frac{1}{p-1} \right) = 0.4574378\ldots
    \end{align}
    Numerically, $B(1,1) = 1.0346534\ldots$
\end{enumerate}
\end{lemma}

\begin{proof}
By the construction in Lemma~\ref{lem:SiegelWalfiszOmega}, we have $\lambda_1(z,g) = A_1(z)\rho(1;z,g)/\Gamma(z)$, where $\rho(1;z,g) = \prod_{p \mid g}(1 - z/p)$ is as in \eqref{eq:rho_def} and
\begin{align}
    A_1(z) := G(1,z) = \prod_p \left(1 - \frac{1}{p}\right)^z \left(1 - \frac{z}{p}\right)^{-1}
\end{align}
denotes the Euler product $G(s,z) := \zeta(s)^{-z}\sum_{n=1}^\infty z^{\Omega(n)} n^{-s}$ evaluated at $s=1$. Note that $A_1(z)$ differs from the function $A_1(z)$ in \cite[Theorem 3]{DupainHallTenenbaum1982} by a factor of $1/\Gamma(z)$.

We begin with the proof of (i). At $z=1$, the product gives $A_1(1) = 1$, and
\begin{align}
    \rho(1;1,g) = \prod_{p \mid g}\left(1 - \frac{1}{p}\right) = \frac{\varphi(g)}{g}.
\end{align}
Thus $b(1,g) = \lambda_1(1,g) = 1 \cdot \frac{\varphi(g)}{g} / 1 = \frac{\varphi(g)}{g}$.

Next, we move on to the proof of (ii). We compute $B(1,g) = \frac{d}{dz}|_{z=1} \lambda_1(z,g)$ using logarithmic differentiation. Since $\lambda_1(z,g) = A_1(z)\rho(1;z,g)/\Gamma(z)$:
\begin{align}
    \frac{d}{dz}\log \lambda_1(z,g) = \frac{A_1'(z)}{A_1(z)} + \frac{\rho'(1;z,g)}{\rho(1;z,g)} - \psi(z)
\end{align}
where $\psi(z) = \Gamma'(z)/\Gamma(z)$ is the digamma function.

At $z=1$, we have $\psi(1) = -\gamma$ and
\begin{align}
    \frac{\rho'(1;z,g)}{\rho(1;z,g)} = -\sum_{p \mid g} \frac{1}{p-z} \quad \Rightarrow \quad \left.\frac{\rho'(1;z,g)}{\rho(1;z,g)}\right|_{z=1} = -\sum_{p \mid g}\frac{1}{p-1}.
\end{align}
For $A_1'(1)/A_1(1) = A_1'(1)$: Taking the logarithmic derivative of $A_1(z)$,
\begin{align}
    \frac{d}{dz}\log A_1(z) = \sum_p \left( \log\left(1-\frac{1}{p}\right) + \frac{1}{p-z} \right).
\end{align}
Once again, at $z=1$ we see that
\begin{align}
    A_1'(1) = \sum_p \left( \log\left(1-\frac{1}{p}\right) + \frac{1}{p-1} \right).
\end{align}
This sum converges since $\log(1-1/p) + 1/(p-1) = O(1/p^2)$ as $p \to \infty$.

Combining, and using $\lambda_1(1,g) = \varphi(g)/g$ leads to
\begin{align}
    B(1,g) = \frac{\varphi(g)}{g}\bigg( A_1'(1) + \gamma - \sum_{p \mid g}\frac{1}{p-1} \bigg) = \frac{\varphi(g)}{g}\bigg( B(1,1) - \sum_{p \mid g}\frac{1}{p-1} \bigg),
\end{align}
which was the last step of the proof.
\end{proof}

Having computed the $j=1$ coefficients, we now show that the higher-order polynomials $P_{j,g}$ collapse to constants, a fact that simplifies the asymptotic analysis of the singular series considerably.

\begin{lemma}[Vanishing of $b(j,g)$ for $j \geq 2$] \label{lem:bj_vanishes}
For all squarefree $g \geq 1$ and all integers $j \geq 2$, we have $b(j,g) = 0$. Consequently, $P_{j,g}(y) = B(j,g)$ is a constant (degree zero polynomial) for all $j \geq 2$.
\end{lemma}

\begin{proof}
The coefficient $b(j,g) = \lambda_j(1,g)$ appears in the Selberg-Delange expansion of Lemma~\ref{lem:SiegelWalfiszOmega}.

At $z = 1$, this becomes
\begin{align}
    \sum_{\substack{n \le x \\ n \equiv h \modu{ug}}} 1 = \frac{x}{\varphi(g)}\left( \lambda_1(1,g) + \frac{\lambda_2(1,g)}{\log x} + \cdots + \frac{\lambda_j(1,g)}{(\log x)^{j-1}} + \cdots \right).
\end{align}

The left-hand side equals $\frac{x}{g} + O(1)$ for $(h,g) = 1$. Since $\lambda_1(1,g) = \varphi(g)/g$ by Lemma~\ref{lem:j1_coefficients}, the leading term matches: $\frac{x}{\varphi(g)} \frac{\varphi(g)}{g} = \frac{x}{g}$.

For consistency with the error $O(1)$, the coefficients of $\frac{x}{\varphi(g)(\log x)^{j-1}}$ must vanish for all $j \geq 1$. Therefore $\lambda_j(1,g) = 0$ for all $j \geq 2$, i.e., $b(j,g) = 0$.
\end{proof}

The preceding lemma shows that $P_{j,g}(y) = B(j,g)$ for $j \ge 2$. We now confirm that these constants are finite, ensuring that the lower-order terms in the expansion of $\mathfrak{f}(q;N,M)/N$ are well-defined.

\begin{lemma}[Finiteness of $B(j,g)$] \label{lem:Bjg_finite}
For all squarefree $g \geq 1$ and all integers $j \geq 1$, the constant $B(j,g) = \frac{d}{dz}|_{z=1} \lambda_j(z,g)$ is finite.
\end{lemma}

\begin{proof}
The Selberg-Delange theory constructs the coefficients $\lambda_j(z,g)$ as functions that are analytic in $z$ in a neighborhood of $z=1$. This analyticity is established in the general theory; see Tenenbaum \cite[Chapter II.5, Theorem 5.2]{tenenbaum}.

Specifically, $\lambda_j(z,g)$ is constructed from the Taylor expansion of the function $H(s,z) = (s-1)^z \sum_{n=1}^\infty z^{\Omega(n)} n^{-s}$ about $s=1$, which is jointly analytic in $(s,z)$ in a neighborhood of $(1,1)$. The coefficients $\lambda_j(z,g)$ inherit analyticity in $z$ from this construction.

Since $\lambda_j(z,g)$ is analytic at $z=1$, its derivative there exists and is finite. Therefore $B(j,g) = \frac{d}{dz}|_{z=1} \lambda_j(z,g)$ is a well-defined finite constant for each $j \geq 1$ and each squarefree $g \geq 1$.
\end{proof}

\subsection{Positivity and asymptotics of the singular series} \label{sec:positivity}
%%%%%%%%%%%%%%%%%%%%%%%%%%%%%%%%%%%%%%%%%%%%%%%%%%%%%%%%%%%%%%%
Armed with the explicit coefficients from the previous subsection, we can now show that the singular series $\mathfrak{S}(N,M)$ is positive for large $N$ and determine its rate of growth. This is the key input for Corollary~\ref{cor:rOmega_asymptotic}.
\begin{theorem}[Positivity and scale of $\mathfrak{S}(N,M)$]  \label{thm:positivity}
Let $M \geq 1$ be a fixed positive integer. There exists $N_0 = N_0(M)$ such that for all $N > N_0$,
\begin{align}
    \mathfrak{S}(N,M) > 0.
\end{align}
Moreover, one has the scale
\begin{align} \label{eq:singular_series_asymp}
    \mathfrak{S}(N,M) = (\log \log N + B(1,1))^3 + O_M(1)
\end{align}
as $N \to \infty$, where $B(1,1) = \gamma + A_1'(1) = 1.034653\ldots$ and the implied constant depends only on $M$.
\end{theorem}

\begin{proof}
We decompose $\mathfrak{S}(N,M) = T_1(N,M) + T_{\geq 2}(N,M)$ where
\begin{align}
    T_1(N,M) = \left( \frac{\mathfrak{f}(1;N,M)}{N} \right)^3
\end{align}
(using $c_1(N) = 1$) and
\begin{align}
    T_{\geq 2}(N,M) = \sum_{q=2}^\infty \left( \frac{\mathfrak{f}(q;N,M)}{N} \right)^3 c_q(N).
\end{align}

The first step is to derive an asymptotic for $T_1(N,M)$.
From the definition of $\mathfrak{f}$,
\begin{align}
    \frac{\mathfrak{f}(1;N,M)}{N} = \sum_{j=1}^M \frac{P_{j,1}(\log \log N)}{(\log N)^{j-1}}.
\end{align}
By Lemma~\ref{lem:j1_coefficients}, the $j=1$ term equals $P_{1,1}(\log\log N) = \log\log N + B(1,1)$. 

By Lemma~\ref{lem:bj_vanishes}, for $j \geq 2$, we have $P_{j,1}(y) = B(j,1)$ (a constant). By Lemma~\ref{lem:Bjg_finite}, each $B(j,1)$ is finite. Therefore, for fixed $M$ we see that
\begin{align}
    \sum_{j=2}^M \frac{P_{j,1}(\log\log N)}{(\log N)^{j-1}} = \sum_{j=2}^M \frac{B(j,1)}{(\log N)^{j-1}} = \frac{B(2,1)}{\log N} + O_M\bigg(\frac{1}{(\log N)^2}\bigg).
\end{align}
Thus
\begin{align} \label{eq:f1_expansion}
    \frac{\mathfrak{f}(1;N,M)}{N} = \log\log N + B(1,1) + \frac{B(2,1)}{\log N} + O_M\bigg(\frac{1}{(\log N)^2}\bigg).
\end{align}

Let $L := \log\log N + B(1,1)$ and $\varepsilon := \frac{B(2,1)}{\log N} + O_M(\frac{1}{(\log N)^2}) = O_M(\frac{1}{\log N})$. Then
\begin{align}
    T_1(N,M) = (L + \varepsilon)^3 = L^3 + 3L^2 \varepsilon + 3L\varepsilon^2 + \varepsilon^3.
\end{align}
Since $L = O(\log\log N)$ and $\varepsilon = O(1/\log N)$ we arrive at $3L^2 \varepsilon = O((\log\log N)^2/\log N) = o(1)$, as well as $3L\varepsilon^2 = O(\frac{\log\log N}{(\log N)^2}) = o(1)$ and $\varepsilon^3 = O(\frac{1}{(\log N)^3}) = o(1)$.
Therefore
\begin{align} \label{eq:T1_asymp}
    T_1(N,M) = (\log\log N + B(1,1))^3 + O_M\bigg(\frac{(\log\log N)^2}{\log N}\bigg).
\end{align}

Next, we examine the cancellation for $q > 1$.
For $q > 1$, the identity $\sum_{g \mid q} \mu(g) = 0$ produces cancellation. The $j=1$ contribution to $\mathfrak{f}(q;N,M)/N$ is
\begin{align}
    \frac{1}{q}\sum_{g \mid q} \frac{g\mu(g)}{\varphi(g)} P_{1,g}\left(\log\log\frac{Ng}{q}\right) = \frac{1}{q}\sum_{g \mid q} \mu(g)\bigg(\log\log\frac{Ng}{q} + B(1,1) - \sum_{p \mid g}\frac{1}{p-1}\bigg).
\end{align}
The $B(1,1)$ terms sum to $B(1,1)\sum_{g \mid q}\mu(g) = 0$.

For the logarithmic terms: writing $\log\log(Ng/q) = \log(\log N + \log(g/q))$, we have for $N$ large and $q \leq N^{1/2}$
\begin{align}
    \log\log\frac{Ng}{q} = \log\log N + \log\left(1 + \frac{\log(g/q)}{\log N}\right) = \log\log N + O\left(\frac{\log q}{\log N}\right).
\end{align}
The $\log\log N$ terms cancel by $\sum_{g \mid q}\mu(g) = 0$, leaving
\begin{align}
    \sum_{g \mid q}\mu(g)\log\log\frac{Ng}{q} = O\left(\frac{d(q)\log q}{\log N}\right).
\end{align}
The remaining sum is bounded: $|\sum_{g \mid q}\mu(g)\sum_{p \mid g}(p-1)^{-1}| \leq d(q)\omega(q)$.

For $j \geq 2$: since $P_{j,g} = B(j,g)$ is constant, these contribute $O_M(d(q)/\log N)$ after summing over $g \mid q$.

Now we need a uniform bound for $|\mathfrak{f}(q;N,M)/N|$ when $q > 1$.
Combining our previous finding on the cancellation for $q>1$, we see that for $2 \leq q \leq N^{1/2}$,
\begin{align}
    \left|\frac{\mathfrak{f}(q;N,M)}{N}\right| \leq \frac{1}{q}\left(|(\Omega*\mu)(q)| + C_M d(q)\omega(q)\right) \leq \frac{C_M' d(q)\omega(q)}{q}
\end{align}
for some constant $C_M'$ depending only on $M$, using $|(\Omega*\mu)(q)| \leq 1$, which holds since $\Omega*\mu$ is supported on prime powers with $(\Omega*\mu)(p^k) = 1$ for all $k \geq 1$ (see Remark~\ref{rem:convolution_structure}).

For $q > N^{1/2}$: by Lemma~\ref{lem:f_bound},
\begin{align}
    \bigg|\frac{\mathfrak{f}(q;N,M)}{N}\bigg| \leq \frac{C_M M\log\log N}{q^{1-\varepsilon}}.
\end{align}

The next order of business consists in bounding $|T_{\geq 2}(N,M)|$.
Using $|c_q(N)| \leq \varphi(q) \leq q$:
for $2 \leq q \leq N^{1/2}$
\begin{align}
    \bigg|\bigg(\frac{\mathfrak{f}(q;N,M)}{N}\bigg)^3 c_q(N)\bigg| \leq \frac{(C_M')^3 d(q)^3\omega(q)^3}{q^2}.
\end{align}
For $q > N^{1/2}$
\begin{align}
    \sum_{q > N^{1/2}} \bigg|\bigg(\frac{\mathfrak{f}(q;N,M)}{N}\bigg)^3 c_q(N)\bigg| \ll_M \frac{(\log\log N)^3}{N^{(1-3\varepsilon)/2}} = o(1).
\end{align}
Since $d(q)^3\omega(q)^3 = O(q^\varepsilon)$ for any $\varepsilon > 0$, the sum $\sum_{q=2}^\infty d(q)^3\omega(q)^3/q^2$ converges. Therefore
\begin{align} \label{eq:T2_bound}
    |T_{\geq 2}(N,M)| \leq C_M''
\end{align}
for some constant $C_M''$ depending only on $M$.

Combining \eqref{eq:T1_asymp} and \eqref{eq:T2_bound} yields
\begin{align}
    \mathfrak{S}(N,M) = (\log\log N + B(1,1))^3 + O_M(1),
\end{align}
where the $O_M(1)$ absorbs both the $O_M((\log\log N)^2/\log N)$ correction from $T_1$ and the bounded contribution from $T_{\geq 2}$.

In particular, $(\log\log N + B(1,1))^3 \to \infty$ as $N \to \infty$, so $\mathfrak{S}(N,M) > 0$ for all $N > N_0(M)$.
\end{proof}

\subsection{Asymptotic for $r_\Omega(N)$}
%%%%%%%%%%%%%%%%%%%%%%%%%%%%%%%%%%%%%%%%%%%%%%%%%%%%%%%%%%%%%%%

We now prove Corollary~\ref{cor:rOmega_asymptotic}.

\begin{proof}[Proof of Corollary~\textnormal{\ref{cor:rOmega_asymptotic}}]
Take $A = 1$ and $\Delta = 1/4$ in Theorem~\ref{thm:partition theorem 1}. Then $M = \lfloor 1 + 16 \cdot 5 \rfloor = 81$, and Theorem 1.3 gives
\begin{align}
    r_\Omega(N) = \frac{N^2}{2}\mathfrak{S}(N,81) + O\left(\frac{N^2(\log\log N)^3}{\log N}\right).
\end{align}
By Theorem~\ref{thm:positivity} with $M = 81$
\begin{align}
    \mathfrak{S}(N,81) = (\log\log N + B(1,1))^3 + O(1).
\end{align}
Substituting we arrive at
\begin{align}
    r_\Omega(N) &= \frac{N^2}{2}\left[(\log\log N + B(1,1))^3 + O(1)\right] + O\bigg(\frac{N^2(\log\log N)^3}{\log N}\bigg) \nonumber \\
    &= \frac{N^2}{2}(\log\log N + B(1,1))^3 + O(N^2),
\end{align}
since the $\frac{N^2}{2} \cdot O(1) = O(N^2)$ term dominates $O(\frac{N^2(\log\log N)^3}{\log N})$.
\end{proof}

\section{Proof of Theorem \ref{thm:omega squared exp sum} and Theorem \ref{thm:omega-quadratic}} \label{sec: proof last 2 theorems}

The two results in this section require fundamentally different techniques. For Theorem~\ref{thm:omega squared exp sum}, the weight $\omega(n)^2$ is not additive, so the methods of Section~\ref{sec:exp_sums_additive} do not apply directly. Instead, we decompose $\omega(n)^2$ into a diagonal term $\omega(n)$ and an off-diagonal sum over pairs of distinct primes, reducing the problem to bounding exponential sums over semiprimes via the results of \cite{DongRoblesZaharescuZeindler}. For Theorem~\ref{thm:omega-quadratic}, the weight $\omega(n)$ is additive but the quadratic phase $\e(\alpha n^2)$ prevents us from using the convolution identity $\omega = \mathbf{1}_{\mathbb{P}} * 1$ effectively; we instead decompose $\omega(n)$ into its mean and fluctuation, handling the former with Weyl's inequality and the latter with the Tur\'{a}n-Kubilius inequality.

\subsection{Proof of Theorem \ref{thm:omega squared exp sum}}
We decompose $\omega(n)^2$ into diagonal and off-diagonal contributions. By opening up the square
\begin{align}
    \omega(n)^2 = \bigg( \sum_{p \mid n} 1 \bigg)^2 = \sum_{p \mid n} \sum_{r \mid n} 1 = \sum_{p \mid n} 1 + 2 \sum_{\substack{p \mid n,\, r \mid n \\ p < r}} 1 = \omega(n) + 2 \sum_{\substack{p < r \\ pr \mid n}} 1,
\end{align}
we arrive at
\begin{align}
    \sum_{n \le X} \omega(n)^2 \e(\alpha n) = \sum_{n \le X} \omega(n) \e(\alpha n) + 2 T(\alpha; X),
\end{align}
where the off-diagonal term is
\begin{align}
    T(\alpha; X) := \sum_{n \le X} \e(\alpha n) \sum_{\substack{p < r \\ pr \mid n}} 1 = \sum_{\substack{p < r \\ pr \le X}} \sum_{m \le X/(pr)} \e(\alpha pr m).
\end{align}
The diagonal is bounded by Remark \ref{remark with Upsilon} with $f = \omega$ and $\Delta = \frac{1}{4}$
\begin{align} \label{eq:omega bound}
\sum_{n \le X} \omega(n) \e(\alpha n) \ll \left( \frac{X}{q^{1/4}} \max\{1, \Upsilon^{1/4}\} + X^{5/6} + X^{3/4} q^{1/4} \right) (\log X)^4.
\end{align} 
It remains to bound $T(\alpha; X)$. Let $U \le X$ be a parameter to be determined. We split $T = T_1 + T_2$ where $T_1$ contains pairs with $pr \le U$ and $T_2$ contains pairs with $pr > U$.

First, we bound $T_1$ (Type I). By the triangle inequality and the geometric series bound,
\begin{align}
    |T_1| \le \sum_{\substack{p < r \\ pr \le U}} \min\left( \frac{X}{pr}, \|\alpha pr\|^{-1} \right) \le \sum_{d \le U} \min\left( \frac{X}{d}, \|\alpha d\|^{-1} \right).
\end{align}
By Lemma 2.2 of \cite{Vaughan1997} together with the extension to $\Upsilon > 0$ as in \cite[Lemma 2.1]{DongRoblesZaharescuZeindler}, we obtain
\begin{align} \label{eq:T1 bound}
    |T_1| \ll \left( \frac{X \max\{1, \Upsilon\}}{q} + U + q \right) (\log X)^2.
\end{align}

Next, we bound $T_2$ (applying the semiprimes bound). For pairs with $pr > U$, we interchange the order of summation. Writing $n = prm$, the constraint $n \le X$ with $pr > U$ becomes $m < X/U$, and for each such $m$ we have $U < pr \le X/m$. Thus
\begin{align}
    T_2 = \sum_{m \le X/U} \sum_{\substack{p < r \\ U < pr \le X/m}} \e(\alpha m  pr).
    \end{align}
Let $S_2(\beta; Y) := \sum_{\substack{p, r \in \Pri \\ pr \le Y}} \e(\beta pr)$ denote the exponential sum over semiprimes as defined in \cite{DongRoblesZaharescuZeindler}. Since $S_2$ counts ordered pairs while our sum has $p < r$, and since the diagonal $p = r$ contributes $O(Y^{1/2})$, we have
\begin{align}
    \sum_{\substack{p < r \\ pr \le Y}} \e(\beta pr) = \frac{1}{2} S_2(\beta; Y) + O(Y^{1/2}).
\end{align}
Therefore
\begin{align}
    T_2 = \frac{1}{2} \sum_{m \le X/U} \left( S_2(\alpha m; X/m) - S_2(\alpha m; U) \right) + O\bigg( \sum_{m \le X/U} (X/m)^{1/2} \bigg).
\end{align}
The error term is $O(X^{1/2} \cdot (X/U)^{1/2}) = O(X/U^{1/2})$, which will be absorbed into the main terms.
    
For each $m$, we apply \cite[Theorem 1.5]{DongRoblesZaharescuZeindler} to bound $S_2(\alpha m; X/m)$. We must first determine the rational approximation for $\beta := m\alpha$. Writing $\alpha = a/q + \lambda$ with $|\lambda| \le \Upsilon/q^2$, we have
\begin{align}
    \beta = m\alpha = \frac{ma}{q} + m\lambda.
\end{align}
Let $g = (m, q)$ and write $m = gm'$, $q = gq'$ with $(m', q') = 1$. Then $ma/q = m'a/q'$ with $(m'a, q') = 1$ (since $(a, q) = 1$ and $(m', q') = 1$). Thus
\begin{align}
    \left| \beta - \frac{m'a}{q'} \right| = |m\lambda| \le \frac{m\Upsilon}{q^2} = \frac{\Upsilon_m}{(q')^2}, \quad \text{where } \Upsilon_m := \frac{m\Upsilon}{g^2} = \frac{m\Upsilon}{(m,q)^2}.
\end{align}
Applying \cite[Theorem 1.5]{DongRoblesZaharescuZeindler} with denominator $q' = q/(m,q)$ and parameter $\Upsilon_m$ gives
\begin{align}
    S_2(\alpha m; X/m) \ll \bigg( \frac{X/m}{(q')^{1/4}} \max\{1, \Upsilon_m^{1/4}\} + \left(\frac{X}{m}\right)^{6/7} + \left( \frac{X}{m} \right)^{3/4} (q')^{1/4} \bigg) (\log X)^3.
\end{align}
Substituting $q' = q/(m,q)$ and simplifying
\begin{align} \label{eq:S2 bound applied}
    S_2(\alpha m; X/m) \ll \bigg( \frac{X (m,q)^{1/4}}{m q^{1/4}} \max\{1, \Upsilon_m^{1/4}\} + \left(\frac{X}{m}\right)^{6/7} + \frac{X^{3/4} q^{1/4}}{m^{3/4} (m,q)^{1/4}} \bigg) (\log X)^3.
\end{align}
    
We first treat the case $\Upsilon \le 1$. Then $\Upsilon_m = m\Upsilon/(m,q)^2 \le m/(m,q)^2 \le m$, so 
\begin{align}
    \max\{1, \Upsilon_m^{1/4}\} \le \max\{1, m^{1/4}\} = m^{1/4} \quad \text{for } m \ge 1.
\end{align}
Inserting this into \eqref{eq:S2 bound applied} and summing over $m \le X/U$

For the first term We require the estimate
\begin{align} \label{eq:gcd sum estimate}
    \sum_{m \le M} \frac{(m,q)^{\frac14}}{m^{\frac34}} = \sum_{g \mid q} g^{\frac14} \sum_{\substack{m \le M \\ (m,q) = g}} \frac{1}{m^{\frac34}} \le \sum_{g \mid q} g^{\frac14} \frac{1}{g^{\frac34}} \sum_{k \le M/g} \frac{1}{k^{\frac34}} \ll \sum_{g \mid q} \frac{(M/g)^{\frac14}}{g^{\frac12}} \ll M^{\frac14} q^\varepsilon
\end{align}
for any $\varepsilon > 0$, where we used $\sum_{g \mid q} g^{-3/4} \ll q^\varepsilon$. Therefore
\begin{align}
    \frac{X}{q^{1/4}} \sum_{m \le X/U} \frac{(m,q)^{1/4}}{m^{3/4}} \ll \frac{X}{q^{1/4}} \frac{(X/U)^{1/4}}{1} q^\varepsilon = \frac{X^{5/4}}{U^{1/4} q^{1/4 - \varepsilon}}.
\end{align}

Next, we focus on the second term for which we have
\begin{align}
    X^{6/7} \sum_{m \le X/U} m^{-6/7} \ll X^{6/7} \left( \frac{X}{U} \right)^{1/7} = \frac{X}{U^{1/7}}.
\end{align}

Lastly, we look at the third term. Similarly to \eqref{eq:gcd sum estimate},
\begin{align}
    X^{3/4} q^{1/4} \sum_{m \le X/U} \frac{1}{m^{3/4} (m,q)^{1/4}} \le X^{3/4} q^{1/4} \sum_{m \le X/U} m^{-3/4} \ll X^{3/4} q^{1/4} \left(\frac{X}{U}\right)^{1/4} = \frac{X q^{1/4}}{U^{1/4}}.
\end{align}

Combining with \eqref{eq:T1 bound}, we have for $\Upsilon \le 1$
\begin{align} \label{eq:T combined bound}
    |T| \ll \bigg( \frac{X}{q} + U + q + \frac{X^{5/4}}{U^{1/4} q^{1/4}} + \frac{X}{U^{1/7}} + \frac{X q^{1/4}}{U^{1/4}} \bigg) (\log X)^3.
\end{align}

The next step is optimization. We apply \cite[Lemma 2.3]{DongRoblesZaharescuZeindler}. Note that $X/q \le X^{5/4}/(U^{1/4} q^{1/4})$ for $U \le X$, and $q \le Xq^{1/4}/U^{1/4}$ for $U \le X$. Thus the dominant terms in \eqref{eq:T combined bound} are
    \begin{align}
        F(U) := U, \quad G_0(U) := \frac{X^{5/4}}{U^{1/4} q^{1/4}}, \quad G_1(U) := \frac{X}{U^{1/7}}, \quad G_2(U) := \frac{X q^{1/4}}{U^{1/4}}.
    \end{align}
These satisfy the hypotheses of \cite[Lemma 2.3]{DongRoblesZaharescuZeindler}: $F$ is increasing, each $G_j$ is decreasing, and the appropriate limits hold. Solving $F(U_j) = G_j(U_j)$
\begin{align}
     U_0^{5/4} = \frac{X^{5/4}}{q^{1/4}} \implies U_0 = \frac{X}{q^{1/5}}, \quad U_1^{8/7} = X \implies U_1 = X^{7/8}, \quad U_2^{5/4} = X q^{1/4} \implies U_2 = X^{4/5} q^{1/5}. \nonumber
\end{align}
By \cite[Lemma 2.3]{DongRoblesZaharescuZeindler}, the minimum of $\max\{F(U), G_0(U), G_1(U), G_2(U)\}$ is achieved when $U = \min \{U_0, U_1, U_2\}$, giving
\begin{align} \label{eq:T bound}
    |T| \ll \left( \frac{X}{q^{1/5}} + X^{7/8} + X^{4/5} q^{1/5} \right) (\log X)^3.
\end{align}
Thus the off-diagonal contribution dominates, and for $\Upsilon \le 1$ we obtain
\begin{align}
    \sum_{n \le X} \omega(n)^2 \e(\alpha n) \ll \left( \frac{X}{q^{1/5}} + X^{7/8} + X^{4/5} q^{1/5} \right) (\log X)^4.
\end{align}

The case $\Upsilon > 1$ is as follows. By Dirichlet's approximation theorem, there exist $a_1 \in \Z$ and $q_1 \in \N$ with $(a_1, q_1) = 1$, $q_1 \le 2q$, and $|\alpha - a_1/q_1| \le 1/(2qq_1) \le 1/q_1^2$. If $a_1/q_1 = a/q$, then $q_1 = q$ and $|\alpha - a/q| \le 1/q^2$, so we may apply the $\Upsilon = 1$ case. Otherwise, $|a_1/q_1 - a/q| \ge 1/(qq_1)$, and
\begin{align}
    \frac{1}{qq_1} \le \left| \frac{a_1}{q_1} - \frac{a}{q} \right| \le |\alpha - a_1/q_1| + |\alpha - a/q| \le \frac{1}{2qq_1} + \frac{\Upsilon}{q^2}.
\end{align}
Thus $1/(2qq_1) \le \Upsilon/q^2$, giving $q^2/(2\Upsilon) \le qq_1 \le 2q^2$, hence $q_1 \ge q/(4\Upsilon)$. Applying the $\Upsilon = 1$ case with denominator $q_1$ yields
\begin{align}
    \sum_{n \le X} \omega(n)^2 \e(\alpha n) &\ll \bigg( \frac{X}{q_1^{1/5}} + X^{7/8} + X^{4/5} q_1^{1/5} \bigg) (\log X)^4 \nonumber \\
        &\ll \bigg( \frac{X \cdot (4\Upsilon)^{1/5}}{q^{1/5}} + X^{7/8} + X^{4/5} (2q)^{1/5} \bigg) (\log X)^4 \nonumber \\
        &\ll \bigg( \frac{X}{q^{1/5}} \Upsilon^{1/5} + X^{7/8} + X^{4/5} q^{1/5} \bigg) (\log X)^4.
\end{align}
This completes the proof.

\subsection{Proof of Theorem \ref{thm:omega-quadratic}}
The convolution approach used for Theorem~\ref{thm:omega squared exp sum} is unavailable here: splitting $\omega(n) = \sum_{p \mid n} 1$ and writing $n = pm$ produces an inner quadratic Weyl sum whose $\Upsilon$-parameter grows with $p$, making the large-prime tail no better than the trivial bound. We circumvent this by writing $\omega(n) = \log\log X + g(n)$ where $g(n) := \omega(n) - \log\log X$ has small variance by the Tur\'{a}n-Kubilius inequality.

\begin{lemma}[Quadratic Weyl bound]\label{lem:weyl}
Let $N\ge 1$ and let $\theta\in\R$. Suppose there exist coprime integers $(b,r)=1$ with $r\ge 1$
such that $|\theta-\frac{b}{r}|\le \frac{1}{r^2}$.
Then
\begin{equation}\label{eq:weyl}
\sum_{n\le N}\e(\theta n^2)
\ll
\frac{N}{\sqrt r}+\sqrt{N\log r}+\sqrt{r\log r}.
\end{equation}
\end{lemma}

This is a standard consequence of Weyl differencing, e.g. \cite[Theorem 4.1]{Vaughan1997} or \cite[Theorem 8.2]{IK04}.

\begin{lemma}[Tur\'{a}n-Kubilius inequality]\label{lem:TK}
There exists an absolute constant $C > 0$ such that for all $X \ge 3$,
\begin{equation}\label{eq:TK}
\sum_{n\le X}(\omega(n)-\log\log X)^2 \le C\, X\log\log X.
\end{equation}
\end{lemma}

This is classical, see Tur\'{a}n \cite{Turan1934} for the original, or \cite[Chapter II.6, Theorem 6.3]{tenenbaum} for a modern treatment. The key input is Mertens' result $\sum_{p \le X} 1/p \sim \log\log X$ with the sum taken over primes.

Let
\[
S:=\sum_{n\le X}\omega(n)\e(\alpha n^2).
\]
We decompose $\omega(n)$ into its mean and fluctuation
\begin{equation}\label{eq:decomposition}
S =(\log\log X)\sum_{n\le X}\e(\alpha n^2) + \sum_{n\le X}g(n)\e(\alpha n^2),
\end{equation}
where $g(n) := \omega(n) - \log\log X$ (note that $g$ depends on the summation range $X$).

For the main term we employ Weyl's bound.
Since $|\alpha - a/q| = |\beta| \le q^{-2}$, Lemma~\ref{lem:weyl} applies directly with $(\theta, b, r) = (\alpha, a, q)$
\begin{equation}\label{eq:main-weyl}
\Bigl|\sum_{n\le X}\e(\alpha n^2)\Bigr|
\ll
\frac{X}{\sqrt q} + \sqrt{X\log q} + \sqrt{q\log q}.
\end{equation}
Multiplying by $\log\log X$ we get
\begin{equation}\label{eq:main-term}
(\log\log X)\Bigl|\sum_{n\le X}\e(\alpha n^2)\Bigr|
\ll
\frac{X\log\log X}{\sqrt q} + (\log\log X)\sqrt{X\log q} + (\log\log X)\sqrt{q\log q}.
\end{equation}

For the residual we make use of Cauchy-Schwarz and Tur\'{a}n-Kubilius. Indeed, for the second sum in \eqref{eq:decomposition}, since $|\e(\alpha n^2)| = 1$, the Cauchy-Schwarz inequality gives
\begin{equation}\label{eq:CS}
\Bigl|\sum_{n\le X}g(n)\e(\alpha n^2)\Bigr|^2
\le
\Bigl(\sum_{n\le X}g(n)^2\Bigr)\Bigl(\sum_{n\le X}1\Bigr)
=
X\sum_{n\le X}g(n)^2.
\end{equation}
By the Tur\'{a}n-Kubilius inequality (Lemma~\ref{lem:TK}) we have
\[
\sum_{n\le X}g(n)^2 = \sum_{n\le X}\bigl(\omega(n)-\log\log X\bigr)^2 \le C\,X\log\log X.
\]
Substituting into \eqref{eq:CS} implies that
\begin{equation}\label{eq:residual}
\Bigl|\sum_{n\le X}g(n)\e(\alpha n^2)\Bigr|
\ll
X\sqrt{\log\log X}.
\end{equation}

Combining \eqref{eq:main-term} and \eqref{eq:residual} via the triangle inequality on \eqref{eq:decomposition}:
\[
|S|
\ll
\frac{X\log\log X}{\sqrt{q}}
+ X\sqrt{\log\log X}
+ (\log\log X)\bigl(\sqrt{X\log q}+\sqrt{q\log q}\,\bigr).
\]
This completes the proof of Theorem~\ref{thm:omega-quadratic}.

%%%%%%%%%%%%%%%%%%%%%%%%%%%%%%%%%%%%%%%%%%%%%%%%%%%%%%%%
\section{Conclusion and future work} \label{sec:conclusion}

A desirable next step after Theorem \ref{thm:exp sum bound class F} would be to enlarge the class $\mathcal{F}_0$.  Consider the more general class $\mathcal{F}_b$, consisting of additive functions $f$ satisfying $f(p) = p^b$.  Notice that 
\begin{align} \label{eq: conclusion beta and B}
    \beta(n) = \sum_{p|n} p \quad \textnormal{and} \quad B(n) = \sum_{p^a\|n} a p.
\end{align} 
are then examples of functions in $\mathcal{F}_1$. It should be possible to modify the methods used in this paper to expand Theorem \ref{thm:exp sum bound class F} to include functions in $\mathcal{F}_b$ for general $b\ge 0$. Moreover, it would also be interesting to refine the $O_M(1)$ error
in \eqref{eq:singular_series_asymp} to achieve a relative error of $O(1/\log N)$. We will address these questions in forthcoming work.

Moreover, rather than considering sums of the form
\begin{align} \label{eq: conclusion sum f}
\sum_{n \le X} f(n) \e(\alpha n)
\end{align} 
where $f$ is an additive function, we could instead consider sums of the form
\begin{align} \label{eq: conclusion z f(n)}
\sum_{n \le X} z^{f(n)} \e(\alpha n)
\end{align} 
where $z \in \C$. This problem is closely related to the generalized divisor function $d_z(n)$ where 
\begin{align} \label{eq: conclusion dz}
d_z(n) = \prod_{p^a || n} \frac{z(z+1)\cdots (z+a-1)}{a!}.
\end{align} 
The Dirichlet series and Euler product associated to $d_z(n)$ are given by
\begin{align} \label{eq: conclusion zetazs}
\zeta^z(s) = \sum_{n=1}^\infty \frac{d_z(n)}{n^s} = \prod_p (1-p^{-s})^{-z}
\end{align} 
for $\real(s)>1$. It would be very interesting to derive a bound for \eqref{eq: conclusion z f(n)}, i.e. find $\mathfrak{f}(X;z)$ such that
\begin{align} \label{eq: conclusion zfn bound}
    \sum_{n \le X} z^{f(n)} \e(\alpha n) \ll \mathfrak{f}(X;z).
\end{align} 
Since $f$ is additive, we see that $z^{f(n)}$ is multiplicative and hence the search for $\mathfrak{f}(X;z)$ could, in principle, be achieved by adapting the techniques from Montgomery and Vaughan \cite{MV77} or Dong \textit{et al.} \cite{DongRoblesZaharescuZeindler}; however tailored approaches are preferred as they lead to tighter bounds. If we recall the work of \cite{SelbergOmega} and \cite{Delange}, then a method of obtaining $\sum_{n \le x} \omega(n)$ involves first finding $\sum_{n \le x} z^{\omega(n)}$ with sufficient uniformity in $z$, then differentiating with respect to $z$, and finally letting $z$ tend to $1$. If one could adapt such techniques, then bounds on \eqref{eq: conclusion sum f} could be achieved in a unified way; however the uniformity on $z$ represents a difficulty.

Moreover, new results on weighted asymptotics could be derived by considering the generating function
\begin{align} \label{eq: conclusion partitions zwf}
\Psi_{w,f}(z) = \sum_{n=1}^\infty \mathfrak{p}_{w,f}(n) w^n = \prod_{n=1}^\infty (1-w^n)^{z^{f(n)}}
\end{align} 
where $w,z \in \C$ and $f$ is an arithmetic, potentially additive, function. To successfully extract $\mathfrak{p}_{w,f}$ through the Hardy-Littlewood circle method, one will need non-trivial bounds of the form \eqref{eq: conclusion zfn bound} to deal with the minor arcs.

Lastly, a second line of future work could be as follows. Assuming the Riemann hypothesis we could integrate $z^{f(n)}$ over the contour in \cite[Figure 2.1]{xmeng2018} that perforates the zeros of $\zeta(s)$ on the critical line and obtain a better understanding of the link the between the zeros of $\zeta(s)$ and the enumeration of partitions weighted by $\omega$ or $\Omega$. In particular, this would allow for an `explicit formula' to arise naturally much like in the prime number theorem. The assumption of the Riemann hypothesis could also be relaxed at the expense of weakening some of the resulting error terms.

\section*{Acknowledgments}
The authors are grateful to Robert Vaughan, G\'{e}rald Tenenbaum, and Owen Sharpe for helpful discussions and feedback on this work.
During the course of this research, AG was partially supported by NSF Grant OIA-2229278 and by Simons Travel Support for Mathematicians Grant No. 946710.

%%%%%%%%%%%%%%%%%%%%%%%%%%%%%%%%%%%%%%%%%%%%%%%%%%%%%%%%
%%%%%%%%%%%%%%%%%%%%%%%%%%%%%%%%%%%%%%%%%%%%%%%%%%%%%%%%
\bibliographystyle{abbrv}
\bibliography{literature_additive}
\end{document}